\documentclass[reqno,twoside,12pt,a4paper]{amsart}
\usepackage{bbm}
\usepackage{amsfonts}
\usepackage{graphicx,amssymb,mathrsfs,amsmath}
\usepackage{amsthm}
\usepackage{hyperref}
\usepackage{soul,color}
\usepackage{geometry}
\usepackage{enumerate}
\usepackage{cancel}
\usepackage{xcolor}
\usepackage[final]{changes}
\definechangesauthor[color=orange]{RevisionRemark}

\allowdisplaybreaks[4]
\newcommand{\N}{\mathbb{N}}
\newcommand{\R}{\mathbb{R}}

\newcommand{\eps}{\mathbb{\varepsilon}}

\newcommand{\dive}{\mathrm{div}}

\newcommand{\pa}{\partial}
\newcommand{\pt}{\partial_t}

\newtheorem{theorem}{Theorem}[section]
\newtheorem{remark}{Remark}[section]

\newtheorem{lemma}{Lemma}[section]
\newtheorem{proposition}{Proposition}[section]

\numberwithin{equation}{section}

\def\div{\mathrm{div}}
\def\v{\nu}
\def\a{\alpha}
\def\r{\varepsilon}

\def\D{\nabla}
\def\fr{\frac{1}{\varepsilon}}
\def\p{\partial_t}
\begin{document}

\title{From bipolar Euler-Poisson system
to unipolar Euler-Poisson system in the perspective of mass  }

\author{Shuai Xi$^1$, Liang Zhao$^{2}$}

\date{}

\maketitle \markboth{S. Xi, L. Zhao}
{From bipolar Euler-Poisson to unipolar one}

\vspace{-3mm}

\begin{center}
{\small $^1$College of Mathematics and Systems Science,\\
Shandong University of Science and Technonlogy, Qingdao 266590, P. R. China\\[2mm]
$^2$School of Mathematical Sciences, Shanghai Jiao Tong University\\
Shanghai 200240, P. R. China}
\end{center}



\vspace{1cm}

\begin{center}
\begin{minipage}{15cm}
\small
\textbf{Abstract.} The main purpose of this paper is to provide an effective procedure to study rigorously the relationship between unipolar and bipolar Euler-Poisson system  in the perspective of mass. Based on the fact that the mass of an electron is far less than that of an ion, we amplify this property by letting $m_e/m_i\rightarrow 0$ and using two different singular limits to illustrate it, which are \added{the} zero-electron mass limit and \added{the} infinity-ion mass limit. We use the method of asymptotic expansion\added{s} to handle the problem and find that the limiting process from bipolar to unipolar system\added{s} is actually the process of decoupling, but not the vanishing of equations of the corresponding \added{the} other particle.
\end{minipage}
\end{center}

\vspace{7mm}

\noindent {\bf Keywords:} {Euler-Poisson system; zero-electron mass limit; infinity-ion mass limit; unipolar; bipolar. }

\vspace{5mm}

\noindent {\bf AMS Subject Classification (2010)~:} 35B40, 35L60, 35Q35

\vspace{5mm}


\section{Introduction}
In the paper, we mainly discuss the fundamental relationship between the unipolar and bipolar system in the perspective of mass based on the famous Euler-Poisson system. \added{As the limiting system of the non-relativistic limit for the Euler-Maxwell system (See \cite{peng2007convergence,yang2010non} and other limiting problems \cite{peng2008convergence,peng2008rigorous,peng2009Asymptotic}), }\replaced{Euler-Poisson system}{which} plays an important role in describing the \replaced{motions}{movement} of charged fluids (ions and electrons) in semi-conductors or plasmas \added{when the effect of the magnetic field is weak}. We consider an un-magnetized plasma consisting of electrons with charge  $ -1 $  and ions with charge  $ +1 $.  More specifically, the scaled Euler-Poisson system in the $d$ dimension\added{al} space $\mathbb{R}^d$ can be described as, with $e$ standing for the electrons and $i$ the ions,
\begin{equation}\label{1epsprt}
\begin{cases}
\partial_tn_{e}+\mathrm{div}\left(n_{e}u_{e}\right)=0,
\\
m_{e}\partial_t\left(n_{e}u_{e}\right)+m_{e}\mathrm{div}\left(n_{e}u_{e}\otimes u_{e}\right)
+\nabla p_e(n_e)= n_e\nabla \phi,\\
\partial_tn_{i}+\mathrm{div}\left(n_{i}u_{i}\right)=0,
\\
m_{i}\partial_t\left(n_{i}u_{i}\right)+m_{i}\mathrm{div}\left(n_{i}u_{i}\otimes u_{i}\right)
+\nabla p_i(n_i)=- n_i\nabla \phi,\\
-\lambda^2\triangle \phi=n_i-n_e,\\
t=0: (n_\nu,u_\nu) = (n_{\nu,0},u_{\nu,0}),\quad \nu=e,i,
\end{cases}
\end{equation}
here for $\nu=$\replaced{$e,i$}{$i,e$},  $ n_\v$ stand for the particle density and $u_\v$ the average velocity \replaced{of}{for} ions and electrons respectively, $ \phi $  is the scaled electric potential. \replaced{All of these are}{These are all} functions of the position  $ x$\added{$=(x_1,\cdots, x_d)$}$ \in \mathbb{R}^d$  and the time  $ t > 0 $.  The pressure functions $ p_\v(n_\v) $
are supposed to be smooth and strictly increasing for all $ n_\v > 0 $.   Usually, they
are of the form
\begin{equation*}
p_\v(n_\v)=a_\v^2n_\v^{\r_\v},\quad \v = e,\ i,
\end{equation*}
\noindent where  $ \r_\v\geq1 $  and  $ a_\v> 0 $  are constants. The fluid is \deleted{called} isothermal if  $ \r_\v=1 $  and adiabatic if  $ \r_\v>1 $.  The parameters $m_\v$ stand for the mass of \added{an} electron\deleted{s} and \added{an} ion\deleted{s} respectively and $\lambda$ \added{$>0$} is the scaled Debye length. For details of the scaling and physical background \added{of the model}, we refer to  \cite{besse2004model,jungel2000hierarchy,chen1984introduction,rishbeth1969introduction} and the reference\added{s} therein. In order to make  $ \phi $ uniquely determined, we add a restriction condition
\begin{equation*}
     \phi(x)\rightarrow 0, \quad  \text{when}  \quad |x|\rightarrow \infty.
\end{equation*}

\noindent \quad  \replaced{Physicists}{Physicians} believe that the \replaced{ions}{electrons} can be regarded as background when studying the equations of \replaced{electrons}{ions} because of the huge mass difference between them. That is to say, \added{the} unipolar model \added{for electrons} was formerly derived from \added{the} bipolar model by assuming that the mass of electrons can be neglected \added{compared to that of ions}. \added{On the other hand, the unipolar model for electrons can also be formally derived based on similar assumptions and simplifications.} However, this \added{kind of simplifications} lack\deleted{s} rigorous proof. To study this, we amplify the \replaced{difference}{relationship} between the mass of ions and electrons by letting $m_e/m_i\rightarrow 0$ and use two different singular limits to illustrate it, which are \added{the} zero-electron mass limit and \added{the} infinity-ion mass limit.  We will prove that the unipolar models are indeed the simplification of the bipolar models. \par

\noindent \quad As is mentioned above,  the study of the limit $m_e/m_i\rightarrow 0$ consists of two natural ways. One is to let $m_i=1$ and  $m_e\rightarrow 0$, which is the known-to-all {\bf zero-electron mass limit}. The limit is based on the assumption that $m_e$ can be ignored when $m_i$ is fixed. Letting $ m_e\rightarrow0 $ and $m_i=1$ in \eqref{1epsprt}, formally we get the system for ions
\begin{equation}
\label{ilimit}
\left\{\aligned
&\partial_tn_{i}+\mathrm{div}\left(n_{i}u_{i}\right)=0,
\\
&\partial_t\left(n_{i}u_{i}\right)+\mathrm{div}\left(n_{i}u_{i}\otimes u_{i}\right)
+\nabla p_i(n_i)=- n_i\nabla \phi,\\
&-\lambda^2\triangle \phi=n_i-n_e,
\endaligned\right.
\end{equation}
and the system for electrons
\begin{equation}\label{elimit}
\left\{\aligned
&\partial_tn_e+\mathrm{div}\left(n_eu_e\right)=0,\\
&\partial_tu_e+\left(u_e\cdot\nabla\right)u_e+\nabla P_e=0,\\
\endaligned\right.
\end{equation}
where $P_e$ is a funtion of $n_e$ and $u_e$. At the same time, we can also obtain  the Maxwell-Boltzmann relationship \cite{0mass},
\[
\nabla p_e(n_e)= n_e\nabla \phi,
\]
which \deleted{together with \eqref{elimit}}is used to replace $n_e$ in \eqref{ilimit}, leading to the solvability of \eqref{ilimit} (see the details in Section 2). We then take back $n_e$ into \eqref{elimit} to solve for $u_e$ and $P_e$, which yields the formal limiting equations for the electrons and success in decoupling. Thus we get the unipolar model of ions \eqref{ilimit} from the bipolar model \eqref{1epsprt}.\par
\noindent \quad Another way is to consider just the opposite, we set $m_e=1$ and $m_i\rightarrow \infty$.  It is based on the fact that $m_i$ turns to infinity when $m_e$ is fixed. We call it \added{the} {\bf infinity-ion mass limit}. We let $m_e=1$ and $m_i\rightarrow +\infty$ in \eqref{1epsprt}, which yields the formal limit system for electrons
\begin{equation}\label{elimit2}
\left\{\aligned
&\partial_tn_{e}+\mathrm{div}\left(n_{e}u_{e}\right)=0,
\\
&\partial_t\left(n_{e}u_{e}\right)+\mathrm{div}\left(n_{e}u_{e}\otimes u_{e}\right)
+\nabla p_i(n_e)= n_e\nabla \phi,\\
&-\lambda^2\triangle \phi=n_i-n_e,
\endaligned\right.
\end{equation}
and the system for ions
\begin{equation}\label{ilimit2}
\begin{cases}
\partial_t n_i+\mathrm{div}(n_iu_i)=0,\\
\partial_t u_i +(u_i\cdot \nabla)u_i=0.\\
\end{cases}
\end{equation}
It is easy to get the \added{local} existence of \added{smooth solutions to} \eqref{ilimit2} by the energy method, then we substitute $n_i$ we have solved in \eqref{ilimit2} into \eqref{elimit2}. The solvability of \eqref{elimit2} is   guaranteed by \added{Lax\cite{La73} and} Kato\cite{Ka75}\deleted{and Majda\cite{Ma84}}. Thus, the decoupling is success\added{ful}. That is to say we get the unipolar model of electrons \eqref{elimit2} from the bipolar model \eqref{1epsprt}. The details of the formal asymptotic analysis can be found in Section 2. \par

The main purpose of this paper is to provide an effective procedure to study rigorously the relationship between unipolar and bipolar systems in the perspective of mass. As to the zero-electron mass limit, many former works have been done \added{for the unipolar Euler-Poisson system for electrons} (see \cite{ali2011zero,ali2010zero,jungel1998zero,xuEP1,Ting0mass}). \added{Due to the complex structure of the bipolar Euler-Poisson system, few results have been obtained by now for the bipolar case. See \cite{TApeng99} for the one-dimensional case in a bounded domain and \cite{lpw} for the quasi-neutral limit. For the unipolar cases, }authors tended to believe that when letting $m_e\rightarrow 0$, the equations of \deleted{the} ions stay the same (see \eqref{ilimit}), so it is rational to ignore the limiting process of the equations of the ions, and put emphasis on the equations of electrons. This is a  misunderstanding. Although the \replaced{equations}{system} for ions \eqref{ilimit} looks the same as \replaced{those}{the equations} of ions in \eqref{1epsprt}, the value of $u_i$ and $n_i$ are different, which actually  are dependent on the parameter $\r\triangleq\sqrt{m_e/m_i}$. Thus, the system for ions \eqref{ilimit} is only invariant in forms. It is improper to just ignore the effect of the ions, and only do the asymptotic analysis to the equations of electrons when considering the two-fluid model. Thus, the limit process from bipolar to unipolar system is actually the process of decoupling, but not the vanishing of equations of the corresponding \added{the} other particle.

The paper is organized as follows. In Section 2, we \deleted{first} introduce some   basic lemmas and give the formal asymptotic analysis as well as the error estimates. The main results of this paper \replaced{are}{is} Theorem \ref{e2.1} and Theorem \ref{i2.1}, which are stated at the end of \added{each subsections in} Section 2. Section 3 and Section 4 are devoted to \added{the} detailed proof\added{s} of Theorem \ref{e2.1} and Theorem \ref{i2.1} in the sense of \added{the} zero-electron mass limit and the infinity-ion mass limit, respectively.

\section{ Preliminaries and main results}
\subsection{ Notations and inequalities.}

In the following, we denote by  $ C $  a generic positive constant independent of $ \r $.  For a multi-index  $ \a=(\a_1,\cdots,\a_d)\in \mathbb{N}^d $  and $ \beta=(\beta_1,\cdots,\beta_d)\in \mathbb{N}^d,\,\beta<\a $  stands for $ \beta\neq\a $  and  $ \beta_j\leq\a_j $  for all  $ j=1,\cdots,d. $  \added{For a multi-index $\a\in\N^d$, we denote
	\[  \pa_x^\alpha=\dfrac{\pa^{|\alpha|}}{\pa x_1^{\alpha_1}\cdots \pa x_d^{\alpha_d}}
	\quad {\text {with}} \quad |\alpha|=\alpha_1+\cdots+\alpha_d.  \]}
\noindent We denote by $ \|\cdot\|_s $,   $ \|\cdot\| $  and  $ \|\cdot\|_\infty $  the norm of the usual Sobolev spaces  $ H^s\left(\mathbb{R}^d\right) $,   $ L^2\left(\mathbb{R}^d\right) $  and  $ L^\infty\left(\mathbb{R}^d\right) $,  respectively. \added{Moreover, we denote for simplicity $H^0=L^2$. }The inner product in $ L^2\left(\mathbb{R}^d\right) $  is denoted by  $ \big\langle\cdot,\,\cdot\big\rangle. $ Throughout the paper, we denote $\nu=e,i$, and
\begin{equation*}
\r=\sqrt{\frac{m_e}{m_i}}.
\end{equation*}

\added{We first give several inequalities that will be used in the later proof.}
\begin{lemma}
\label{lem2.2} (Moser-type calculus inequalities, see \cite{klainerman1981singular} and \cite{Ma84} ).
Let  $ s \geq 1 $  be an integer. Suppose  $ u\in H^s\left(\mathbb{R}^d\right),\,\D u\in L^\infty\left(\mathbb{R}^d\right) $  and  $ v\in H^{s-1}\left(\mathbb{R}^d\right)\cap L^\infty\left(\mathbb{R}^d\right) $. Then for all  $ \a\in\mathbb{N}^d $  with  $ 1\leq|\a|\leq s $  and all smooth function  $ f $,  we have  $ \partial^\a_x(uv)-u\partial^\a_xv\in L^2\left(\mathbb{R}^d\right) $,   $ \partial^\a_xf(u)\in L^2\left(\mathbb{R}^d\right) $  and
\begin{eqnarray*}
 \|\partial^\a_x(uv)-u\partial^\a_xv\|&\leq& C_s\left(\|\D u\|_{\infty}\|\D^{|\a|-1}v\|+\|\D^{|\a|}u\|\,\| v\|_{\infty}\right), \\
   \|\partial^\a_xf(u)\|&\leq& C_\infty\left(\|\D u\|_{\infty}+1\right)^{|\a|-1}\|\D^{|\a|}u\|,
\end{eqnarray*}
where the constant  $ C_\infty>0 $  depends on  $ \| u\|_{\infty} $  and  $ s $,  and  $ C_s>0 $  \added{is a generic constant which }depends only on   $ s $.
Moreover, if  $ s > \frac{d}{2}+1 $,  then the embedding  $ u\in H^s\left(\mathbb{R}^d\right)\hookrightarrow W^{1,\infty}\left(\mathbb{R}^d\right) $  is continuous and \added{thus} we have
 $$  \|\partial^\a_x(uv)-u\partial^\a_xv\|\leq C_s\|\D u\|_{s-1}\|v\|_{s-1}. $$
\end{lemma}

\begin{lemma}
	\label{lem2.1}
	Let $s > \frac{d}{2}+2$ be an integer \deleted{and $d\leq 3$}. For all $\a\in\mathbb{N}^d$ with $2\leq|\a|\leq s$,
	if $u\in H^s\left(\R^d\right)$ and $v\in H^{s}\left(\R^d\right)$,
	then
	\begin{eqnarray*}
		\|\partial^\a_x(uv)-u\partial^\a_xv-\sum_{{\substack{1\leq i\leq d\\\alpha_i\neq 0}}  }\alpha_i\pa_{x_i} u\pa^{\alpha^i}_x v\| &\leq& C_s\left(\|\D ^2 u\|_{\infty}\|\D^{|\a|-2}v\|+\|\D^{|\a|}u\|\| v\|_{\infty}\right)\nonumber\\
		&\leq& C_s\|\D^2 u\|_{s-2}\|v\|_{s-2},
	\end{eqnarray*}
	where $\alpha^i\in\N^d$ is a multi-index and $\pa_{x_i}\pa^{\alpha^i}_x=\pa^{\alpha}_x $. Term $\displaystyle\sum_{{\substack{1\leq i\leq d\\\alpha_i\neq 0}}  }\alpha_i\pa_{x_i} u\pa^{\alpha^i}_x v$ denotes all the terms related to the first order derivatives of $u$ by using the Leibniz Formula\deleted{s}.
\end{lemma}
\begin{lemma}
\label{lem2.3}
 For any smooth \added{vector-valued} function  $ u $ \deleted{$:\,\mathbb{R}^d\rightarrow\mathbb{R}^d $}  and \added{scalar function} $ \Phi$ \deleted{$:\,\mathbb{R}^d\rightarrow\mathbb{R} $},  we have
 $$  |\left<u\triangle\Phi,\nabla\Phi\right>|\leq C\|\D u\|_{\infty}\|\D \Phi\|^2, $$
where the constant  $ C>0 $  is independent of  $ u $  and  $ \Phi $.
\end{lemma}

In order to \replaced{simplify}{simply} the later proof, we introduce the enthalpy function, defined as
\begin{equation*}
h_\nu'(n)=\frac{p_\nu'(n)}{n} \mbox{\quad and \quad}h_\nu(1)=0.
\end{equation*}
Then for \replaced{$n_\nu^{m_e,m_i}>0$}{$n_\nu>0$}, system \eqref{1epsprt} can be rewritten into
\begin{equation}\label{1}
\begin{cases}

\pt n_e^{m_e,m_i} +\div(n_e^{m_e,m_i}u_e^{m_e,m_i})=0,\\
\pt u_e^{m_e,m_i}+(u_e^{m_e,m_i}\cdot\nabla)u_e^{m_e,m_i}+\displaystyle\frac{\nabla h_e(n_e^{m_e,m_i})}{m_e}=\displaystyle\frac{\nabla \phi^{m_e,m_i}}{m_e}, \\
\pt n_i^{m_e,m_i}+\div(n_i^{m_e,m_i}u_i^{m_e,m_i})=0,\\
\pt u_i^{m_e,m_i}+(u_i^{m_e,m_i}\cdot\nabla)u_i^{m_e,m_i}+\displaystyle\frac{\nabla h_i(n_i^{m_e,m_i})}{m_i}=-\displaystyle\frac{\nabla \phi^{m_e,m_i}}{m_i},\\
-\lambda^2\Delta \phi^{m_e,m_i}=n_i^{m_e,m_i}-n_e^{m_e,m_i},\\
t=0:  (n_\nu^{m_e,m_i},u_\nu^{m_e,m_i}) = (n_{\nu,0}^{m_e,m_i},u_{\nu,0}^{m_e,m_i}),\quad \nu=e,i,
\end{cases}
\end{equation}
\added{where the initial data of $\phi^{m_e,m_i}$, i.e.$\phi_0^{m_e,m_i}$, is defined by
\[
-\lambda^2 \Delta \phi_0^{m_e,m_i}=n_{i,0}^{m_e,m_i}-n_{e,0}^{m_e,m_i}.
\]
}

The next result concerns the local existence of smooth solutions which can be
easily obtained by employing the theory of \added{Lax\cite{La73} and }Kato\cite{Ka75} for the
symmetrizable hyperbolic system.

\begin{proposition}
\label{2.4}
Let  $ s > \dfrac{d}{2}+1 $  be an integer and $ \left(n_{\nu,0}^{m_e,m_i},\,u_{\nu,0}^{m_e,m_i}\right)\in H^s\left(\mathbb{R}^d\right) $  with  $ n_{\nu,0}^{m_e,m_i}\geq 2\underline{n} $  for some given constant  $ \underline{n}>0 $,  independent of  $m_e$ \added{and} $m_i$.  Then there exists  $ T^{m_e,m_i}_1>0 $  such that the Cauchy problem \eqref{1} has a unique smooth solution  $ \left(n^{m_e,m_i}_\nu,\,u^{m_e,m_i}_\nu,\,\phi^{m_e,m_i}\right) $  defined \replaced{on the}{in} time interval $ [0,T^{m_e,m_i}_1] $, satisfying  $ n^{m_e,m_i} \geq \underline{n} $  and
\begin{eqnarray*}
\left(n^{m_e,m_i}_\nu,\,u^{m_e,m_i}_\nu\right)&\in &C\left(\left[0,\,T^{m_e,m_i}_1\right];\,H^s\left(\mathbb{R}^d\right)\right)\cap C^1\left(\left[0,\,T^{m_e,m_i}_1\right];\,H^{s-1}\left(\mathbb{R}^d\right)\right), \nonumber\\
 \phi^{m_e,m_i}&\in& C\left(\left[0,\,T^{m_e,m_i}_1\right];\,H^{s+1}\left(\mathbb{R}^d\right)\right)\cap C^1\left(\left[0,\,T^{m_e,m_i}_1\right];\,H^{s}\left(\mathbb{R}^d\right)\right).
\end{eqnarray*}
\end{proposition}

\subsection{\bf \large  Asymptotic analysis for zero-electron mass limit ($\r\rightarrow 0, m_i=1$)}
\subsubsection{Formal expansion}

As to the zero-electron mass limit, \added{by} setting $m_i=1$, we look for an approximation of solution  $ \left(n_\nu^{\r,1},\,u_\nu^{\r,1},\,\phi^{\r,1}\right) $  to \eqref{1} in the form of power series \added{with respect to the small parameter $\eps$}.  \added{In this subsection, we denote the integer $s> \frac{d}{2}+2$.} Assume that the initial data of  $ \left(n_\nu^{\r,1},\,u_\nu^{\r,1},\,\phi^{\r,1}\right) $  admit an asymptotic expansion with respect to  $ \r $, \deleted{for $\nu=e,i$,}
\begin{equation*}\label{4.2}
\left(n_{\nu,0}^{\r,1},\,u_{\nu,0}^{\r,1},\,\phi_0^{\r,1}\right)(x)
=\sum_{j\geq0}\r^{2j}\left(\bar{n}_\nu^{e,j},\,\bar{u}_\nu^{e,j},\,\bar{\phi}^{e,j}\right)(x), \added{\quad \nu=e,i,}
\end{equation*}
where  $ \left(\bar{n}_\nu^{e,j},\,\bar{u}_\nu^{e,j},\,\bar{\phi}^{e,j}\right)_{j\geq0} $
are sufficiently smooth\replaced{.}{,} \added{We further assume} \deleted{and} the following ansatz:
\begin{equation}
\label{4.3}
\left(n_\nu^{\r,1},\,u_\nu^{\r,1},\,\phi^{\r,1}\right)(t,\,x)
=\sum_{j\geq0}\r^{2j}\left(n_\nu^{e,j},\,u_\nu^{e,j},\,\phi^{e,j}\right)(t,\,x).
\end{equation}
In what follows, we use a formal expansion formula, \added{which is obtained by the Taylor's formula}
\begin{equation*}
h_\nu\left(\sum_{j\geq0}\r^{2j}n_\nu^{e,j}\right)=h_\nu\left(n_\nu^{e,0}\right)+h_\nu'\left(n_\nu^{e,0}\right)\sum_{j\geq1}\r^{2j}n_\nu^{e,j}+\sum_{j\geq2}\r^{2j}h_\nu^{e,j-1}\left(\left(n_\nu^{e,k}\right)_{k\leq j-1}\right),
\end{equation*}
where $ \{h_\nu^{e,j}\}_{j\geq 1} $ are smooth function\added{s} depending only on $h_\nu$  and \replaced{$ \left(n_{\nu}^{e,k}\right)_{k\leq j} $}{$ \left(n^k\right)_{k\leq j} $}. \added{We also denote $h_\nu^{e,0}=0$.} For simplicity, from now on, we denote $ \left(n_\nu^{\r,1},\,u_\nu^{\r,1},\,\phi^{\r,1}\right) $\replaced{,}{and} $ \left(n_\nu^{e,j},\,u_\nu^{e,j},\,\phi^{e,j}\right)_{ j\geq 0} $, \added{ $\left(\bar{n}_\nu^{e,j},\,\bar{u}_\nu^{e,j},\,\bar{\phi}^{e,j}\right)_{j\geq 0}$ and  $ \{h_\nu^{e,j}\}_{j\geq 0} $} by
$ \left(n_\nu^\r,\,u_\nu^\r,\,\phi^\r\right) $ \replaced{,}{and}  $ \left(n_\nu^j,\,u_\nu^j,\,\phi^j\right)_{ j\geq 0} $, \added{$\left(\bar{n}_\nu^{j},\,\bar{u}_\nu^{j},\,\bar{\phi}^{j}\right)_{j\geq 0}$ and $ \{h_\nu^{j}\}_{j\geq 0} $ respectively} in the section of zero-electron mass limit. Substituting the expansions \eqref{4.3} into system \eqref{1} \added{and comparing the coefficients before each order of $\eps$, }we obtain
\vspace{5mm}

\noindent (1) The leading profiles  $ \left(n_\nu^0,\,u_\nu^0,\,\phi^0\right) $  satisfy the following system 
\begin{equation}\label{0E-M}
\begin{cases}
\pt n_e^0+\div(n_e^0u_e^0)=0,\\
\pt u_e^0+(u_e^0\cdot\nabla)u_e^0+\nabla P_e^0=0,\\
\pt n_i^0+\div(n_i^0u_i^0)=0,\\
\pt u_i^0+(u_i^0\cdot\nabla)u_i^0+\nabla(h_i(n_i^0)+\phi^0)=0,\\
-\lambda^2\Delta \phi^0=n_i^0-n_e^0,
\end{cases}
\end{equation}
where
\begin{equation}\label{P0}
P_e^0=h_e^\prime(n_e^0)n_e^1-\phi^1,
\end{equation}
with the initial data
\begin{equation}
\label{4.6}
   \left(n_\nu^0,u_\nu^0\right)(0,x) = \left(\bar{n}_\nu^0,\bar{u}_\nu^0\right)(x),
\quad x\in\mathbb{R}^d.
\end{equation}
\replaced{Noticing that}{Notice} the \added{coefficients before the }$\r^{-2}$ term \added{imply}
\begin{equation}\label{Boltzmann}
\nabla h_e(n_e^0)-\nabla \phi^0=0,
\end{equation}
we deduce that \added{up to a constant} $n_e^0=h_e^{-1}(\phi^0)$. Thus the equations for ions and the Poisson equation in \eqref{0E-M} are actually the following unipolar Euler-Poisson system for ions
\begin{equation}\label{e2.9}
\begin{cases}
\pt n_i^0+\div(n_i^0u_i^0)=0,\\
\pt(u_i^0)+(u_i^0\cdot\nabla)u_i^0+\nabla(h_i(n_i^0)+\phi^0)=0,\\
-\lambda^2\Delta \phi^0=n_i^0-h_e^{-1}(\phi^0).
\end{cases}
\end{equation}
The solvability of the Poisson equation can be found in \cite{lgp}, in which $\phi^0$ is expressed as a function of $n_i^0$ \added{and can be viewed as a zeroth order term}. Thus the first two equations are \added{symmetrizable} hyperbolic\deleted{system}, \replaced{to}{of} which the unique local smooth solution exists due to the famous work of \added{Lax\cite{La73} and} Kato\cite{Ka75}\deleted{ and Majda\cite{Ma84}}. \replaced{By this time}{At the same time} $n_e^0$ is also known since it is a function of $\phi^0$, and $ (u_e^0,\,P_e^0) $  satisfy the following
incompressible Euler equations:
\begin{equation}\label{e2.10}
\left\{\aligned
&\div\left(n_e^0u_e^0\right)=-\partial_tn_e^0,\\
&\partial_tu_e^0+\left(u_e^0\cdot\D\right)u_e^0+\D P_e^0=0,
\endaligned
\right.
\end{equation}
\added{to which the local existence of smooth solutions is ensured by the theories of Lax\cite{La73} and Kato\cite{Ka75}}.Then Cauchy problem \eqref{0E-M} with \eqref{4.6}
has then been solved, \added{and} $ P_e^0 $  is used to solve \deleted{the}  $ \left(n_e^1,\,\phi^1\right) $  by \eqref{P0}.

\vspace{5mm}
\noindent(2) \replaced{For}{At} the  $ \varepsilon^2$  term, we find
\begin{equation}\label{1E-M}
\left\{\aligned
&\partial_tn_e^1+\div\left(n_e^0u_e^1+n_e^1u_e^0\right)=0,\\
&\partial_tu_e^1+\left(u_e^0\cdot\D\right)u_e^1+\left(u_e^1\cdot\D\right)u_e^0+\D P_e^1=0,\\
&\partial_tn_i^1+\div\left(n_i^0u_i^1+n_i^1u_i^0\right)=0,\\
&\partial_tu_i^1+\left(u_i^0\cdot\D\right)u_i^1+\left(u_i^1\cdot\D\right)u_i^0+\D (h_i^\prime(n_i^0)n_i^1+\phi^1)=0,\\
&-\lambda^2\Delta \phi^1=n_i^1-n_e^1,
\endaligned
\right.
\end{equation}
where
\begin{equation*}
P_e^1=h_e^\prime(n_e^0)n_e^2+h_e^1(n_e^1)-\phi^2,
\end{equation*}
with the  initial data
\begin{equation}
\label{p4.6}
   \left(n_\nu^1,u_\nu^1\right)(0,x) = \left(\bar{n}_\nu^1,\bar{u}_\nu^1\right)(x),
\quad x\in\mathbb{R}^d.
\end{equation}
Since $P_e^0$ \added{and $n_e^0$ are} \deleted{is} known, \added{by \eqref{P0},} we substitute $\phi^1$ into the Poisson equation in \eqref{1E-M} with
\begin{equation*}
n_e^1=\frac{P_e^0+\phi^1}{h_e^\prime(n_e^0)}.
\end{equation*}
\deleted{which implies $n_e^1$ is a function of $n_i^1$.}Now the equations for ions in \eqref{1E-M} turn out to be the following linear system
\begin{equation*}
\begin{cases}
\partial_tn_i^1+\div\left(n_i^0u_i^1+n_i^1u_i^0\right)=0,\\
\partial_tu_i^1+\left(u_i^0\cdot\D\right)u_i^1+\left(u_i^1\cdot\D\right)u_i^0+\D (h_i^\prime(n_i^0)n_i^1+\phi^1)=0,\\
-\added{\lambda^2}\Delta \phi^1=n_i^1-\displaystyle\frac{P_e^0+\phi^1}{h_e^\prime(n_e^0)},
\end{cases}
\end{equation*}
\replaced{to}{for} which we  can get the unique solution $n_i^1,u_i^1$ and $\phi^1$, and thus $n_e^1$. Also, $u_e^1$ and $P_e^1$ satisfy the following
\begin{equation*}
\begin{cases}
\partial_tn_e^1+\div\left(n_e^0u_e^1+n_e^1u_e^0\right)=0,\\
\partial_tu_e^1+\left(u_e^0\cdot\D\right)u_e^1+\left(u_e^1\cdot\D\right)u_e^0+\D P_e^1=0,\\
\end{cases}
\end{equation*}
in which $P_e^1$ is used to solve $(u_e^2,\phi^2)$. \par

\vspace{5mm}
\noindent (3) For  $ j\geq2 $, in general the profiles  $ \left(n_\nu^j,\,u_\nu^j,\,\phi^j\right) $
are obtained by induction. Assume that
 $ \left(n_\nu^k,\,u_\nu^k,\,\phi^k\right)_{0\leq k\leq j-1} $  are smooth and already
determined in previous steps. Then  $ \left(n_\nu^j,\,u_\nu^j,\,\phi^j\right) $
satisfy the linear system
\begin{equation}
\label{jE-M}
\left\{\aligned
&\partial_tn_e^j+\div\left(n_e^0u_e^j+n_e^ju_e^0\right)=-\sum_{k=1}^{j-1}\div\left(n_e^ku_e^{j-k}\right),\\
&\partial_tu_e^j+\left(u_e^0\cdot\D\right)u_e^j+\left(u_e^j\cdot\D\right)u_e^0+\D P_e^j=-\displaystyle\sum_{k=1}^{j-1}(u_e^k\cdot\nabla)u_e^{j-k},\\
&\partial_tn_i^j+\div\left(n_i^0u_i^j+n_i^ju_i^0\right)=-\sum_{k=1}^{j-1}\div\left(n_i^ku_i^{j-k}\right),\\
&\partial_tu_i^j+\left(u_i^0\cdot\D\right)u_i^j+\left(u_i^j\cdot\D\right)u_i^0+\D P_i^j=-\displaystyle\sum_{k=1}^{j-1}(u_i^k\cdot\nabla)u_i^{j-k},\\
&-\lambda^2\triangle\phi^j=n_i^j-n_e^j,
\endaligned\right.
\end{equation}
where
\begin{equation*}
\begin{cases}
P_e^j=h_e'(n_e^0)n_e^{j+1}+h_e^j\left((n_e^k)_{k\leq j}\right)-\phi^{j+1},\\
P_i^j=h_i'(n_i^0)n_i^j+h_i^{j-1}\left((n_i^k)_{k\leq j-1}\right)+\phi^{j},
\end{cases}
\end{equation*}
with the initial data
\begin{equation}
\label{j2.5}
\left(n_\nu^j,\,u_\nu^j\right)(0,\,x)
=\left(\bar{n}_\nu^j,\,\bar{u}_\nu^j\right)(x),\quad \v=e,\,i.
\end{equation}
Generally, we can get  $ \left(n_i^j,\ u_i^j,\ \phi^j\right) $  from  $ P_e^{j-1} $  and the third to fifth equations in \eqref{jE-M}, \added{and thus $n_e^j$. }\deleted{and}Then we can obtain $ (u_e^j,\,P_e^j) $  from the first two equations in \eqref{jE-M}. \added{We then have the following Proposition.}
\begin{proposition}
Assume that the initial data $ \left(\bar{n}_\nu^j,\,\bar{u}_\nu^j,\,\bar{\phi}^j\right)_{j\geq 0} $  are sufficiently smooth with  $ \bar{n}_\nu^0 > 0 $  in  $ \mathbb{R}^d $.  Then there exist \deleted{the} unique smooth profiles  $ \left(n_\nu^j,\,u_\nu^j,\,\phi^j\right)_{j\geq 0} $,  \added{which are obtained as} solutions \replaced{to}{of the}  problems \eqref{0E-M} with \eqref{4.6},  \eqref{1E-M} with \eqref{p4.6} and \eqref{jE-M} with \eqref{j2.5} \replaced{on}{in} the time interval  $ [0,T_1^e] $ {\textbf{\added{with $T_1^e$ independent of $\eps$.}}} In other words, there exists a unique asymptotic expansion up to any order of the form \eqref{4.3}.
\end{proposition}

\subsubsection{Error estimates   and main result}
Let  \replaced{$m\geq1$}{$ m\in \N $}  be a fixed integer and  $ \left(n_\nu^\r,\,u_\nu^\r,\,\phi^\r\right) $  be the exact solution to problem \eqref{1} (with $m_i=1$) defined \replaced{on the}{in} time interval $ \left[0,\,T_1^{\r,1 }\right] $.  We denote by $ \left(n^{\added{e,}m}_{\nu,\r},\,u^{\added{e,}m}_{\nu,\r},\,\phi^{\added{e,}m}_{\r}\right) $  the approximate solution of order $ m $  defined \replaced{on}{in}  $ \left[0,\,T_1^e \right] $  by
\begin{equation*}
\left(n^{\added{e,}m}_{\nu,\r},\,u^{\added{e,}m}_{\nu,\r},\,\phi^{\added{e,}m}_{\r}\right)=\sum_{j=0}^{m}\r^{2j}\left(n_\nu^j,\,u_\nu^j,\,\phi^j\right),
\end{equation*}
where  $ \left(n_\nu^j,\,u_\nu^j,\,\phi^j\right)_{0\leq j\leq m} $  are constructed in the previous subsection. The proof of the convergence of the asymptotic expansion \eqref{4.3} is to establish the limit
 $$       \left(n_\nu^\r,\,u_\nu^\r,\,\phi^\r\right)
-\left(n^{\added{e,}m}_{\nu,\r},\,u^{\added{e,}m}_{\nu,\r},\,\phi^{\added{e,}m}_{\r}\right) \longrightarrow 0,   $$
and \added{obtain} its convergence rate as  $ \r\rightarrow0 $  \replaced{on}{in} a time interval independent of
 $ \r $,  when the convergence holds at  $ t=0 $. \added{For convenience, we denote $\left(n^{\added{e,}m}_{\nu,\r},\,u^{\added{e,}m}_{\nu,\r},\,\phi^{\added{e,}m}_{\r}\right)$ by $\left(n^{m}_{\nu,\r},\,u^{m}_{\nu,\r},\,\phi^{m}_{\r}\right)$ in the section of zero-electron mass limit.}

For  $ \v= e,\,i $,  \added{we} define the \replaced{remaining terms}{remainders} $ \left({R}_{n_\nu}^{\varepsilon,1,m},\,R_{u_\nu}^{\varepsilon,1,m}\right) $  by
\begin{equation}
\label{E-M2}
\left\{\aligned
&\partial_tn_{e,\r}^{m}+\div\left(n_{e,\r}^{m}u_{e,\r}^{m}\right)=R_{n_e}^{\varepsilon,1,m},
 \\
&\partial_tu_{e,\r}^{m}+\left(u_{e,\r}^{m}\cdot\D\right)u_{e,\r}^{m}
+\frac{1}{\r^2}\D \left(h_e(n_{e,\r}^{m})-\phi^m_{\r}\right)=\added{\dfrac{1}{\eps^2}}R_{u_e}^{\varepsilon,1,m}, \\
&\partial_tn_{i,\r}^{m}+\div\left(n_{i,\r}^{m}u_{i,\r}^{m}\right)=R_{n_i}^{\varepsilon,1,m},
 \\
&\partial_tu_{i,\r}^{m}+\left(u_{i,\r}^{m}\cdot\D\right)u_{i,\r}^{m}
+\D \left(h_i(n_{i,\r}^{m})+\phi^m_{\r}\right)=R_{u_i}^{\varepsilon,1,m}, \\
&-\lambda^2\triangle\phi_{\r}^{m}
=n_{i,\r}^{m}-n_{e,\r}^m.
\endaligned\right.
\end{equation}
It is clear that the convergence rate depends strongly on the order of the \replaced{remaining terms}{remainders} with respect to  $ \r $.  Since the profiles $ \left(n_\nu^j,\,u_\nu^j,\,\phi^j\right)_{0\leq j\leq m} $  are sufficiently smooth, we have

\begin{proposition}\label{error}
   If \eqref{0E-M}, \eqref{1E-M} and \eqref{jE-M} hold, then we can find $\widetilde{R}_{n_{\nu}}^{\varepsilon,1,m}$, such that
\begin{equation*}
\div{\widetilde{R}_{n_{\nu}}^{\varepsilon,1,m}}= R_{n_{\nu}}^{\varepsilon,1,m},
\end{equation*}
and for all integers  \replaced{$ m \geq 1 $,}{$ m \geq 0 $}\deleted{ and $ s \geq 0 $}  the \replaced{remaining terms}{remainders} satisfy
\begin{equation}
\label{2.11}
\sup_{0\leq t\leq T_1^e}\left\|\left(\widetilde{R}_{n_{\nu}}^{\varepsilon,1,m},\added{R_{n_{\nu}}^{\varepsilon,1,m},}\,R_{u_\nu}^{\r,1,m}\right)(t)\right\|_s\leq C_m \r^{2m\added{+2}},
\end{equation}
where  $ C_m > 0 $  is a constant independent of  $ \r $.
\end{proposition}
\noindent \textbf{Proof.} By the definition of $R_{n_\v}^{\varepsilon,1,m}$ in \eqref{E-M2}, we have
\begin{eqnarray*}
R_{n_{\nu}}^{\varepsilon,1,m}&=&\partial_tn_{{\nu},\varepsilon}^m
  +\mathrm{div}\left(n_{{\nu},\varepsilon}^{m}u_{{\nu},\varepsilon}^m\right)
  \nonumber\\
&=&\sum_{j=0}^{m}\varepsilon^{2j}\partial_tn^j_{\nu}
+\mathrm{div}\left(\left(\sum_{j=0}^{m}\varepsilon^{2j}n^j_{\nu}\right)\left(\sum_{j=0}^{m}
\varepsilon^{2j}u^j_{\nu}\right)\right)\nonumber\\
&=&\partial_tn_\nu^0+\mathrm{div}\left(n_\nu^0u_\nu^0\right)
+\sum_{j=1}^{m}\varepsilon^{2j}\left(\partial_tn^j_{{\nu}}+\sum_{k=0}^{j}
\mathrm{div}\left(n^k_{{\nu}}u^{j-k}_{{\nu}}\right)\right) \nonumber\\
&\quad& +\sum_{j=1}^{m}\varepsilon^{2j+2m}\left(\sum_{k=j}^{m}
\mathrm{div}\left(n^k_{{\nu}}u^{m+j-k}_{{\nu}}\right)\right)
\nonumber\\
&=&\mathrm{div}\left(\sum_{j=1}^{m}\varepsilon^{2j+2m}\left(\sum_{k=j}^{m}
n^k_{{\nu}}u^{m+j-k}_{{\nu}}\right)\right)\overset{\added{\Delta}}{=}\mathrm{div}\left( \widetilde{R}_{n_{\nu}}^{\varepsilon,1,m}\right),
\end{eqnarray*}
then $$    \sup_{0\leq t\leq T_1^e}\|\widetilde{R}_{n_{\nu}}^{\varepsilon,1,m}(t)\|_s \leq C_m\varepsilon^{2m+2}.   $$
By the definition of $R_{u_e}^{\varepsilon,1,m}$, we have
\begin{eqnarray*}
\added{\dfrac{1}{\eps^2}}R_{u_e}^{\r,1,m}&=&\partial_tu_{e,\r}^m+\left(u_{e,\r}^{m}\cdot\D\right)u_{e,\r}^{m}
+\frac{1}{\r^2}\D \left(h_e(n_{e,\r}^{m})-\phi^m_{\r}\right)\nonumber\\
&=&\sum_{j=0}^{m}\r^{2j}\partial_tu^j_e+\left(\left(\sum_{j=0}^m\r^{2j}u_e^j\right)\cdot\nabla\right)\sum_{j=0}^m\r^{2j}u_e^j\nonumber\\
&\quad&+\frac{1}{\r^2}\nabla\left(h_e\left(n_e^0\right)+h_e'\left(n_e^0\right)\sum_{j\geq1}\r^{2j}n_e^j+\sum_{j\geq2}\r^{2j}h_e^{j-1}\left(\left(n_e^k\right)_{k\leq j-1}\right)\right)-\frac{1}{\r^2}\sum_{j=0}^m\r^{2j}\D\phi^j\nonumber\\
&=&\dfrac{1}{\r^2}\nabla (h_e(n_e^0)-\phi^0)+\left(\pt u_e^0+(u_e^0\cdot\nabla)u_e^0+\nabla \left(h_e^\prime(n_e^0)n_e^1-\phi^1\right)\right)\nonumber\\
&\quad&+\sum_{j=1}^{m-1}\r^{2j}\left(\pt u_e^j+\sum_{k=0}^{j}(u_e^k\cdot\nabla)u_e^{j-k}+\nabla\left(h_e^\prime(n_e^0)n_e^{j+1}+h_e^{j}\left(\left(n_e^k\right)_{k\leq j}\right)-\phi^{j+1}\right)\right)\nonumber\\
&\quad&+O(\r^{2m}),
\end{eqnarray*}
and by the definition of $R_{u_i}^{\varepsilon,1,m}$, we have
\begin{eqnarray*}
\qquad R_{u_i}^{\r,1,m}&=&\partial_tu_{i,\r}^m+\left(u_{i,\r}^{m}\cdot\D\right)u_{i,\r}^{m}
+\D \left(h_i(n_{i,\r}^{m})+\phi^m_{\r}\right)\nonumber\\
&=&\sum_{j=0}^{m}\r^{2j}\partial_tu^j_i+\left(\left(\sum_{j=0}^m\r^{2j}u_i^j\right)\cdot\nabla\right)\sum_{j=0}^m\r^{2j}u_i^j\nonumber\\
&\quad&+\nabla\left(h_i\left(n_i^0\right)+h_i'\left(n_i^0\right)\sum_{j\geq1}\r^{2j}n_i^j+\sum_{j\geq2}\r^{2j}h_i^{j-1}\left(\left(n_i^k\right)_{k\leq j-1}\right)\right)+\sum_{j=0}^m\r^{2j}\D\phi^j\nonumber\\
&=&\pt u_i^0+(u_i^0\cdot\nabla)u_i^0+h_i(n_i^0)+\nabla \phi^0\nonumber\\
&\quad&+\r^2\left(\pt u_i^0+(u_i^0\cdot\nabla)u_i^1+(u_i^1\cdot\nabla)u_i^0+\nabla \left(h_i^\prime(n_i^0)n_i^1+\phi^1\right)\right)\nonumber\\
&\quad&+\sum_{j=2}^m\r^{2j}\left(\pt u_i^j+\sum_{k=0}^{j}(u_i^k\cdot\nabla)u_i^{j-k}+\nabla\left(h_i^\prime(n_i^0)n_i^j+h_i^{j-1}\left(\left(n_i^k\right)_{k\leq j-1}\right)+\phi^j\right)\right)\nonumber\\
&\quad&+O(\r^{2m+2}).
\end{eqnarray*}
Hence, \added{combining} \eqref{0E-M}, \eqref{1E-M}, \eqref{jE-M} and the Maxwell-Boltzmann relationship \eqref{Boltzmann}, \replaced{implies}{imply} \eqref{2.11}.      \hfill$\square$

The main result for the zero-electron mass limit is the following convergence result, of which the proof
will be given in  Section 3.
\begin{theorem} [The zero-electron mass limit]
\label{e2.1}Under the conditions of Proposition \ref{error},
   let  $ s > \frac{d}{2}+2 $  and  \replaced{$ m\geq1 $}{$ m\in\N $}  be integers. Assume
\begin{eqnarray}
\label{e2.12}
&\,&\added{\|n_{e,0}^{\r,1}-n_{e,\r}^{e,m}(0,\,\cdot)\|_s+\|n_{i,0}^{\r,1}-n_{i,\r}^{e,m}(0,\,\cdot)\|_s}\nonumber\\
&\added{+}&\added{{\r}\|u_{e,0}^{\r,1}-u_{e,\r}^{e,m}(0,\,\cdot)\|_s+\|u_{i,0}^{\r,1}-u_{i,\r}^{e,m}(0,\,\cdot)\|_s\leq C_1\eps^{2m\added{+1}},}
\end{eqnarray}
where  $ C_1 > 0 $  is a constant independent of  $ \r $.  Then, for the isothermal fluid, there exists a constant  $ C_2 > 0 $,  which depends on  $ T_1^e$  but is independent of  $ \r $,  such that as  $ \r\rightarrow0 $  we have  $ T_1^{\r,1}\geq T_1^e $   and for all integers  $ 2m> s $,  the solution  $ \left(n_{\added{\nu}}^{\r,1},\,u_{\added{\nu}}^{\r,1},\,\phi^{\r,1}\right)$  to the problem \eqref{1} satisfies
\begin{eqnarray*}\label{e2.13}
&\,&\added{\sup_{0\leq t\leq T_1^e}\left(\|n_e^{\r,1}(t)-n_{e,\r}^{e,m}(t)\|_s^2+\|n_i^{\r,1}(t)-n_{i,\r}^{e,m}(t)\|_s^2+\|\nabla (\phi^{\eps,1}(t)-\phi_\eps^{e,m}(t)\|_s^2\right)}\nonumber\\
&\,&\added{+\sup_{0\leq t\leq T_1^e}\left({\r}\|u_e^{\r,1}(t)-u_{e,\r}^{e,m}(t)\|_s^2+\|u_i^{\r,1}(t)-u_{i,\r}^{e,m}(t)\|_s^2\right)\leq  C_2\r^{2(2m\added{+1}-s)}.}
\end{eqnarray*}
That is to say, the zero-electron mass limit  $ \r\rightarrow 0 $  of the bipolar  Euler-Poisson system  \eqref{1} is the unipolar  Euler-Poisson equations for ions \eqref{e2.9} and the incompressible Euler equations \eqref{e2.10}.
\end{theorem}
\subsection{\bf \large  Asymptotic analysis for infinity-ion mass limit ($\r\rightarrow 0, m_e=1$)}
\subsubsection{Formal expansion}
As to the infinity-ion mass limit, \added{by} setting $m_e=1$, we look for an approximation of solution  $ \left(n_\nu^{1,\frac{1}{\r}},\,u_\nu^{1,\frac{1}{\r}},\,\phi^{1,\frac{1}{\r}}\right) $  to \eqref{1} in the form of power series.  \added{In this subsection, we denote the integer $s> \frac{d}{2}+1$.} Assume that the initial data of  $ \left(n_\nu^{1,\frac{1}{\r}},\,u_\nu^{1,\frac{1}{\r}},\,\phi^{1,\frac{1}{\r}}\right) $  admit an asymptotic expansion with respect to  $ \r $, \deleted{for $\nu=e,i$,}
\begin{equation*}\label{5.2}
\left(n_{\nu,0}^{1,\frac{1}{\r}},\,u_{\nu,0}^{1,\frac{1}{\r}},\,\phi_0^{1,\frac{1}{\r}}\right)(x)
=\sum_{j\geq0}\r^{2j}\left(\bar{n}_\nu^{i,j},\,\bar{u}_\nu^{i,j},\,\bar{\phi}^{i,j}\right)(x), \added{\quad \nu=e,i, }
\end{equation*}
where  $ \left(\bar{n}_\nu^{i,j},\,\bar{u}_\nu^{i,j},\,\bar{\phi}^{i,j}\right)_{j\geq0} $
are sufficiently smooth. \added{We further assume} \deleted{and} the following ansatz:
\begin{equation}
\label{5.3}
\left(n_\nu^{1,\frac{1}{\r}},\,u_\nu^{1,\frac{1}{\r}},\,\phi^{1,\frac{1}{\r}}\right)(t,\,x)
=\sum_{j\geq0}\r^{2j}\left(n_\nu^{i,j},\,u_\nu^{i,j},\,\phi^{i,j}\right)(t,\,x).
\end{equation}
In what follows, we use the formal expansion \replaced{which is obtained by the Taylor's formula}{defined by}
\begin{equation*}
h_\nu\left(\sum_{j\geq0}\r^{2i}n_\nu^{i,j}\right)=h_\nu\left(n_\nu^{i,0}\right)+h_\nu'\left(n_\nu^{i,0}\right)\sum_{j\geq1}\r^{2j}n_\nu^{i,j}+\sum_{j\geq2}\r^{2j}h_\nu^{i,j-1}\left(\left(n_\nu^{i,k}\right)_{k\leq j-1}\right),
\end{equation*}
where $\{ h_\nu^{i,j}\}_{j\geq 1} $  are smooth functions depending only on \added{$h_\nu$ and} $ \left(n_{\added{\nu}}^{\added{i,} k}\right)_{k\leq j} $. \added{We also formally denote $h_\nu^{i,-1}=h_\nu^{i,0}=0$.}
For simplicity, from now on, we denote $ \left(n_\nu^{1,\frac{1}{\r}},\,u_\nu^{1,\frac{1}{\r}},\,\phi^{1,\frac{1}{\r}}\right) $\replaced{,}{and} $\left(n_\nu^{i,j},\,u_\nu^{i,j},\,\phi^{i,j}\right)_{ j\geq 0} $, \added{$\left(\bar{n}_\nu^{i,j},\,\bar{u}_\nu^{i,j},\,\bar{\phi}^{i,j}\right)_{j\geq 0}$ and $\{ h_\nu^{i,j}\}_{j\geq -1} $} by $ \left(n_\nu^\r,\,u_\nu^\r,\,\phi^\r\right) $\replaced{,}{and} $ \left(n_\nu^j,\,u_\nu^j,\,\phi^j\right)_{ j\geq 0} $, \added{\\\noindent$\left(\bar{n}_\nu^{j},\,\bar{u}_\nu^{j},\,\bar{\phi}^{j}\right)_{j\geq 0}$ and $\{ h_\nu^{j}\}_{j\geq -1} $} in the \added{section of} infinity-ion mass limit. Substituting the expansions \eqref{5.3} into system \eqref{1} \added{and comparing the coefficients before each order of $\eps$,} we obtain

\noindent(1) The leading profiles  $ \left(n_\nu^0,\,u_\nu^0,\,\phi^0\right) $  satisfy the following system 
\begin{equation}\label{0I-M}
\begin{cases}
\pt n_i^0+\div(n_i^0u_i^0)=0,\\
\pt u_i^0+(u_i^0\cdot\nabla)u_i^0=0,\\
\pt n_e^0+\div(n_e^0u_e^0)=0,\\
\pt u_e^0+(u_e^0\cdot\nabla)u_e^0+\nabla h_e(n_e^0)= \nabla \phi^0,\\
-\lambda^2\Delta \phi^0=n_i^0-n_e^0,
\end{cases}
\end{equation}
with the initial data
\begin{equation}
\label{5.6}
   \left(n_\nu^0,u_\nu^0\right)(0,x) = \left(\bar{n}_\nu^0,\bar{u}_\nu^0\right)(x),
\quad x\in\mathbb{R}^d.
\end{equation}
Through energy method, it is easy to get the unique \added{local smooth} solution $(n_i^0, u_i^0)$ \replaced{to}{of} the following system
\begin{equation}\label{i2.29}
\begin{cases}
\pt n_i^0+\div(n_i^0u_i^0)=0,\\
\pt u_i^0+(u_i^0\cdot\nabla)u_i^0=0.
\end{cases}
\end{equation}
Since $n_i^0$ is known, we can see that the third to fifth equation in \eqref{0I-M} is actually the decoupled unipolar Euler-Poisson system for electrons,
\begin{equation}\label{i2.30}
\begin{cases}
\pt n_e^0+\div(n_e^0u_e^0)=0,\\
\pt u_e^0+(u_e^0\cdot\nabla)u_e^0+\nabla h_e(n_e^0)= \nabla \phi^0,\\
-\lambda^2\Delta \phi^0=n_i^0(t,x)-n_e^0,
\end{cases}
\end{equation}
and thus $(n_e^0,u_e^0,\phi^0)$ is known due to \added{Lax\cite{La73} and} Kato\cite{Ka75}\deleted{and Majda\cite{Ma84}}.
\vspace{5mm}

\noindent(2) In general, for  $ j\geq1 $,  the profiles  $ \left(n_\nu^j,\,u_\nu^j,\,\phi^j\right) $
are obtained by induction. Assume that
 $ \left(n_\nu^k,\,u_\nu^k,\,\phi^k\right)_{0\leq k\leq j-1} $  are smooth and already
determined in previous steps. Then  $ \left(n_\nu^j,\,u_\nu^j,\,\phi^j\right) $
satisfy the linear system
\begin{equation}
\label{jI-M}
\left\{\aligned
&\partial_tn_i^j+\div\left(n_i^0u_i^j+n_i^ju_i^0\right)=-\sum_{k=1}^{j-1}\div\left(n_i^ku_i^{j-k}\right),\\
&\partial_tu_i^j+\left(u_i^0\cdot\D\right)u_i^j+\left(u_i^j\cdot\D\right)u_i^0+\D P_i^j=-\displaystyle\sum_{k=1}^{j-1}(u_i^k\cdot\nabla)u_i^{j-k},\\
&\partial_tn_e^j+\div\left(n_e^0u_e^j+n_e^ju_e^0\right)=-\sum_{k=1}^{j-1}\div\left(n_e^ku_e^{j-k}\right),\\
&\partial_tu_e^j+\left(u_e^0\cdot\D\right)u_e^j+\left(u_e^j\cdot\D\right)u_e^0+\D P_e^j=-\displaystyle\sum_{k=1}^{j-1}(u_e^k\cdot\nabla)u_e^{j-k},\\
&-\lambda^2\triangle\phi^j=n_i^j-n_e^j,
\endaligned\right.
\end{equation}
where
\begin{equation*}
\begin{cases}
P_e^j=h_{\added{e}}'(n_e^0)n_e^{j}+h_{\added{e}}^{j-1}\left((n_e^k)_{k\leq j-1}\right)-\phi^{j},\\
P_i^j=h_{\added{i}}'(n_i^0)n_i^{j-1}+h_{\added{i}}^{j-2}\left((n_i^k)_{k\leq j-2}\right)+\phi^{j-1},
\end{cases}
\end{equation*}
with the initial data
\begin{equation}
\label{j3.5}
\left(n_\nu^j,\,u_\nu^j\right)(0,\,x)
=\left(\bar{n}_\nu^j,\,\bar{u}_\nu^j\right)(x),\quad \v=e,\,i.
\end{equation}
\added{We mention here that formally $h_\nu^{-1}=h_\nu^{0}=0$.} Generally, we can get  $ \left(n_i^j,\ u_i^j\right) $  from  the first two equations in \eqref{jI-M}, and then insert \replaced{$n_i^j$}{$u_i^j$} into the \added{Poisson equation to solve the }third to fifth equations in \eqref{jI-M} \replaced{and thus getting}{to} $ (n_e^j,\,u_e^j,\ \phi^j) $. \added{We then have the following proposition.}
\begin{proposition}
Assume that the initial data $ \left(\bar{n}_\nu^j,\,\bar{u}_\nu^j,\,\bar{\phi}^j\right)_{j\geq 0} $  are sufficiently smooth with  $ \bar{n}_\nu^0 > 0 $  in  $ \mathbb{R}^d $.  Then there exist \deleted{the} unique smooth profiles  $ \left(n_\nu^j,\,u_\nu^j,\,\phi^j\right)_{j\geq 0} $,  \added{which are obtained as} solutions \replaced{to}{of the} problems \eqref{0I-M} with \eqref{5.6} and \eqref{jI-M} with \eqref{j3.5} \replaced{on}{in} the time interval  $ [0,T_1^i] $ {\textbf{\added{with $T_1^i$ independent of $\eps$}}}. That is to say there exists a unique asymptotic expansion up to any order of the form \eqref{5.3}.
\end{proposition}

\subsubsection{Error estimates   and main result}
Let  \replaced{$m\geq1$}{$ m\in \N $}  be a fixed integer and  $ \left(n_\nu^\r,\,u_\nu^\r,\,\phi^\r\right) $  be the exact solution to problem \eqref{1} defined \replaced{on the}{in} time interval $ \left[0,\,T_1^{1,\fr} \right] $.  We denote by $ \left(n^{\added{i,}m}_{\nu,\r},\,u^{\added{i,}m}_{\nu,\r},\,\phi^{\added{i,}m}_{\r}\right) $  the approximate solution of order $ m $  defined \replaced{on}{in}  $ \left[0,\,T_1^i\right] $  by
\begin{equation*}
\left(n^{\added{i,}m}_{\nu,\r},\,u^{\added{i,}m}_{\nu,\r},\,\phi^{\added{i,}m}_{\r}\right)=\sum_{j=0}^{m}\r^{2j}\left(n_\nu^j,\,u_\nu^j,\,\phi^j\right),
\end{equation*}
where  $ \left(n_\nu^j,\,u_\nu^j,\,\phi^j\right)_{0\leq j\leq m} $  are constructed in the previous subsection. The convergence of the asymptotic expansion \eqref{5.3} is to establish the limit
 $$       \left(n_\nu^\r,\,u_\nu^\r,\,\phi^\r\right)
-\left(n^{\added{i,}m}_{\nu,\r},\,u^{\added{i,}m}_{\nu,\r},\,\phi^{\added{i,}m}_{\r}\right) \longrightarrow 0,  $$
and \added{obtain} its convergence rate as  $ \r\rightarrow0 $  \replaced{on}{in} a time interval independent of
 $ \r $,  when the convergence holds at  $ t=0 $. \added{For convenience, we denote $\left(n^{\added{i,}m}_{\nu,\r},\,u^{\added{i,}m}_{\nu,\r},\,\phi^{\added{i,}m}_{\r}\right)$ by $\left(n^{m}_{\nu,\r},\,u^{m}_{\nu,\r},\,\phi^{m}_{\r}\right)$ in the section of infinity-ion mass limit.}

For  $ \v= e,\,i $,  define the \replaced{remaining terms}{remainders} $ \left(R_{n_\nu}^{1,\fr,m},\,R_{u_\nu}^{1,\fr,m}\right) $  by
\begin{equation*}
\left\{\aligned
&\partial_tn_{e,\r}^{m}+\div\left(n_{e,\r}^{m}u_{e,\r}^{m}\right)=R_{n_e}^{1,\fr,m},
 \\
&\partial_tu_{e,\r}^{m}+\left(u_{e,\r}^{m}\cdot\D\right)u_{e,\r}^{m}
+\D \left(h_{\added{e}}(n_{e,\r}^{m})-\phi^m_{\r}\right)=R_{u_e}^{1,\fr,m}, \\
&\partial_tn_{i,\r}^{m}+\div\left(n_{i,\r}^{m}u_{i,\r}^{m}\right)=R_{n_i}^{1,\fr,m},
 \\
&\partial_tu_{i,\r}^{m}+\left(u_{i,\r}^{m}\cdot\D\right)u_{i,\eps}^{m}
+\r^2\D \left(h_{\added{i}}(n_{i,\r}^{m})+\phi^m_{\r}\right)=R_{u_i}^{1,\fr,m}, \\
&-\lambda^2\triangle\phi_{\r}^{m}
=n_{i,\r}^{m}-n_{e,\r}^m.
\endaligned\right.
\end{equation*}
It is clear that the convergence rate depends strongly on the order of the \replaced{remaining terms}{remainders} with respect to  $ \r $.  Since the profiles $ \left(n_\nu^j,\,u_\nu^j,\,\phi^j\right)_{0\leq j\leq m} $  are sufficiently smooth, a straightforward computation gives the following result.

\begin{proposition}\label{error1}
   Let \eqref{0I-M} and \eqref{jI-M} hold. \added{There exist a smooth function $\widetilde{R}_{n_\nu}^{1,\fr,m}$ such that\[\dive \widetilde{R}_{n_\nu}^{1,\fr,m}={R}_{n_\nu}^{1,\fr,m},
   	\] and} for all integers  \replaced{$m\geq1$}{$ m \geq 0 $}\deleted{and
 $ s \geq 0 $},  the \replaced{remaining terms}{remainders} satisfy
\begin{equation*}
\label{2.11im}
\sup_{0\leq t\leq T_1^{\added{i}}}\big\|\big(\widetilde{R}_{n_\nu}^{1,\fr,m},\,\added{{R}_{n_\nu}^{1,\fr,m},\,}R_{u_\nu}^{1,\fr,m}\big)(t)\big\|_s\leq C_m\r^{2m+2},
\end{equation*}
where  $ C_m > 0 $  is a constant independent of  $ \r $.
\end{proposition}
\noindent The proof is similar to \added{that of} Proposition \ref{error}, we omit it here. The main result for the infinity-ion mass limit is the following convergence result, of which the proof will be given in  Section 4.
\begin{theorem}[The infinity-ion mass limit]
\label{i2.1}Under the conditions of Proposition \ref{error1},
   Let  $ s > \frac{d}{2}+1 $  and  \replaced{$m\geq1$} {$ m \in \N $}  be integers. Assume
\begin{eqnarray*}
\label{i2.12}
&\quad&\added{\left\|n_{e,0}^{1,\fr}-n_{e,\r}^{i,m}(0,\,\cdot)\right\|_s+\left\|n_{i,0}^{1,\fr}-n_{i,\r}^{i,m}(0,\,\cdot)\right\|_s}\nonumber\\
&\quad&\added{+\left\|\frac{1}{\r}\left(u_{i,0}^{1,\fr}-u_{i,\r}^{i,m}(0,\,\cdot)\right)\right\|_s+\left\|u_{e,0}^{1,\fr}-u_{e,\r}^{i,m}(0,\,\cdot)\right\|_s\leq C_1 \eps^{2m+2},}
\end{eqnarray*}
where  $ C_1 > 0 $  is a constant independent of  $ \r $.  Then  there exists a constant  $ C_2 > 0 $,  which depends on  $ T_1^i $  but is independent of  $ \r $,  such that as  $ \r\rightarrow0 $  we have  $ T_1^{1,\fr}\geq T_1^i $   and for all integers  \replaced{$m\geq1$}{$ m\geq 0 $},  the solution  $ \left(n^{1,\fr},\,u^{1,\fr},\,\phi^{1,\fr}\right)$ to the problem \eqref{1} satisfies
\begin{eqnarray*}
\label{i2.13}
&\quad&\added{\sup_{0\leq t\leq T_1^i}\left(\left\|n_e^{1,\fr}(t)-n_{e,\r}^{i,m}(t)\right\|_s^2+\left\|n_i^{1,\fr}(t)-n_{i,\r}^{i,m}(t)\right\|_s^2+\left\|\D\left(\phi^{1,\fr}(t)-\phi^{i,m}_\r(t)\right)\right\|_s^2\right)}\nonumber\\
&\quad&\added{+\sup_{0\leq t\leq T_1^i}\left(\left\|u_e^{1,\fr}(t)-u_{e,\r}^{i,m}(t)\right\|_s^2+\dfrac{1}{\eps^2}\left\|u_i^{1,\fr}(t)-u_{i,\r}^{i,m}(t)\right\|_s^2\right)\leq C_2\eps^{4m\added{+2}}.}
\end{eqnarray*}
That is to say, the infinity-ion mass limit  $ \r\rightarrow 0 $  of the bipolar  Euler-Poisson system  \eqref{1} is the unipolar  Euler-Poisson equations for electrons \eqref{i2.30} and the  equations \eqref{i2.29}.
\end{theorem}
\begin{remark}We mention the difference of \added{the} condition\added{s} needed for \added{the} zero-electron mass limit and the infinity-ion mass limit. In Theorem \ref{e2.1}, we require the fluid to be isothermal and the integer $s>\frac{d}{2}+2$, which is like the situation in \cite{0mass}. And in Theorem \ref{i2.1}, the \deleted{isothermal} condition \added{that the fluid should be isothermal} is not needed, and $s>\frac{d}{2}+1$.
\end{remark}

\section{Proof of Theorem \ref{e2.1}}
\subsection{Energy estimates.} In this section, {\added{we study the zero-electron mass limit of the bipolar Euler-Poisson equations. We assume the fluid to be isothermal, i.e., $p_\nu(n_\nu^{\eps,1})=a_\nu^2 n_\nu^{\eps,1}$. In this section,}}  we continue to use $ \left(n_\nu^\r,\,u_\nu^\r,\,\phi^\r\right) $ and $ \left(n_\nu^j,\,u_\nu^j,\,\phi^j\right)_{ j\geq 0} $ to substitute $ \left(n_\nu^{\r,1},\,u_\nu^{\r,1},\,\phi^{\r,1}\right) $ and $ \left(n_\nu^{e,j},\,u_\nu^{e,j},\,\phi^{e,j}\right)_{ j\geq 0} $. \added{All the corresponding notations are in accordance with those in Subsection 2.2}. The exact solution $\left(n_\v^{\r},\,u_\v^{\r},\,\phi^{\r}\right)$
is defined \replaced{on the}{in} time interval $\left[0,\,T_1^{\r,1} \right]$ and the approximate
solution $\left(n^m_{\v,\r},\,u^m_{\v,\r},\,\phi^m_\r\right)$  \replaced{on the}{in} time
interval $\left[0,\,T_1^e \right]$, with  \replaced{$T_1^e$ independent of $\eps$.}{$T_1^{\r,1}>0$ and $T_1^e>0$.} Let
$$      T_2^{\r,1}=\min\left(T_1^{\r,1},\,T_1^e\right)>0,   $$
then the exact solution and the approximate solution are both defined \replaced{on the}{in} time interval $\left[0,\,T_2^{\r,1} \right]$. \replaced{On}{In} this time interval, we denote
\begin{equation*}\label{3.1}
\left(N_\v^{\r},\,U_\v^{\r},\Phi^{\r}\right)
\triangleq\left(n_\v^\r-n_{\v,\r}^{m},\,u_\v^\r-u_{\v,\r}^{m},\phi^\r-\phi^m_\r\right),\quad \v=e,\,i.
\end{equation*}
\deleted{For simplicity, we denote
$\left(N_\v^{\r,1},\,U_\v^{\r,1},\Phi^{\r,1},R^{\r,1,m}\right)$ by $\left(N_\v^{\r},\,U_\v^{\r},\Phi^{\r},R^{\r,m}\right)$ in this section.} It is easy to check that \deleted{ the variable}
$\left(N_\v^\r,\,U_\v^\r\right)$ satisfy
\begin{equation}
\label{E-M3}
\left\{\aligned
&\p N_{e}^{\r}+\big(U_{e}^{\r}+u_{e,\r}^{m}\big)\D N_{e}^{\r}
+\big(N_{e}^{\r}+n_{e,\r}^{m}\big)\div U_{e}^{\r}
=-\big(N_{e}^{\r}\div u_{e,\r}^{m}+U_{e}^{\r}\D n_{e,\r}^{m}\big)-R_{n_e}^{\r,m}, \\
&\r\p U_{e}^{\r}+\r\big(\big(U_{e}^{\r}+u_{e,\r}^{m}\big)\cdot\D\big) U_e^\r
+\frac{1}{\r}h'_{e}(N_e^\r+n_{e,\r}^{m})\D N_e^\r\\
&\hspace{1.5cm}=-\r\big(U_{e}^{\r}\cdot\D\big) u_{e,\r}^m
-\frac{1}{\r}\big(h'_{e}(N_e^\r+n_{e,\r}^m)-h'_{e}(n_{e,\r}^m)\big)\D n_{e,\r}^m  +\frac{1}{\r}\nabla\Phi^\r-\frac{1}{\r}R_{u_e}^{\r,m},\\
&\p N_{i}^{\r}+\big(U_{i}^{\r}+u_{i,\r}^{m}\big)\D N_{i}^{\r}
+\big(N_{i}^{\r}+n_{i,\r}^{m}\big)\div U_{i}^{\r}
=-\big(N_{i}^{\r}\div u_{i,\r}^{m}+U_{i}^{\r}\D n_{i,\r}^{m}\big)-R_{n_i}^{\r,m}, \\
&\p U_{i}^{\r}+\big(\big(U_{i}^{\r}+u_{i,\r}^{m}\big)\cdot\D\big) U_i^\r
+h'_{i}(N_i^\r+n_{i,\r}^{m})\D N_i^\r\\
&\hspace{1.5cm}=-\big(U_{i}^{\r}\cdot\D\big) u_{i,\r}^m
-\big(h'_{i}(N_i^\r+n_{i,\r}^m)-h'_{i}(n_{i,\r}^m)\big)\D n_{i,\r}^m -\nabla\Phi^\r-R_{u_i}^{\r,m},\\
&\big.\big(N_\v^\r,\,U_\v^\r\big)\big|_{t=0} =\big(n_{\v,0}^\r-n_{\v,\r}^m(0,\,\cdot),
\,u_{\v,0}^\r-u_{\v,\r}^m(0,\,\cdot)\big),
\endaligned
\right.
\end{equation}
coupled with the Poisson equation for  $ \Phi^\r $
\begin{equation*}
\label{F3.2}
    -\lambda^2\Delta \Phi^\r=N_{i}^{\r}-N_{e}^{\r}, \qquad \lim_{|x|\rightarrow +\infty}\Phi^\r(x) = 0.
\end{equation*}
For simplicity, we let $\lambda=1$. Set
$$  W_{ e}^{\r}=\left(\begin{array}{c}
     N_{e}^{\r} \\
      {\r}U_{e}^{\r} \\
   \end{array}
 \right),\quad
W_{i}^{\r}=\left(\begin{array}{c}
     N_{i}^{\r} \\
     U_{i}^{\r} \\
   \end{array}
 \right),
$$
\begin{eqnarray*}
H_{e,\r}^{1}&=&
\left(
   \begin{array}{c}
   N_{e}^{\r}\div u_{e,\r}^{m}+U_{e}^{\r}\D n_{e,\r}^{m} \\
\r\left(U_{e}^{\r}\cdot\D\right) u_{e,\r}^{m}+
\frac{1}{\r}\left(h'_{e}(N_e^\r+n_{e,\r}^{m})-h'_{e}(n_{e,\r}^{m})\right)\D n_{e,\r}^{m}\\
   \end{array}
 \right),
\nonumber\\
\\
H_{i,\r}^{1}&=&
\left(
   \begin{array}{c}
   N_{i}^{\r}\div u_{i,\r}^{m}+U_{i}^{\r}\D n_{i,\r}^{m} \\
\left(U_{i}^{\r}\cdot\D\right) u_{i,\r}^{m}+\left(h'_{i}(N_i^\r+n_{i,\r}^{m})-h'_{i}(n_{i,\r}^{m})\right)\D n_{i,\r}^{m} \\
   \end{array}
 \right),
\end{eqnarray*}
$$
H_{e,\r}^{2}=
\left(
   \begin{array}{c}
0 \\
\frac{1}{\r}\nabla\Phi^\r\\
   \end{array}
 \right),
\quad
H_{i,\r}^{2}=
\left(
   \begin{array}{c}
0 \\
-\nabla\Phi^\r\\
   \end{array}
 \right),
$$$$
R_{e}^{\r}=
\left(
   \begin{array}{c}
    R_{n_e}^{\r} \\
\frac{1}{\r}R_{u_e}^{\r} \\
   \end{array}
 \right) ,
 R_{i}^{\r}=
\left(
   \begin{array}{c}
    R_{n_i}^{\r} \\
R_{u_i}^{\r} \\
   \end{array}
 \right), $$
and for \replaced{$1\leq j\leq d$}{$j=1,2,3$} \added{and $u_\nu^\eps=(u_{\nu,1}^\eps, \cdots, u_{\nu,d}^\eps)$,}
\begin{eqnarray*}
A_e^j\left(n_{e}^{\r},\,u_{e}^{\r}\right)&=&
\left(
  \begin{array}{cc}
u_{e,j}^{\r} & \frac{1}{\r}n_{e}^{\r}e_j^{\top} \\
\frac{1}{\r}h'_e\left(n_{e}^{\r}\right)e_j& u_{e,j}^{\r}\mbox{\bf I}_d \\
  \end{array}
\right),\nonumber\\
A_i^j\left(n_{i}^{\r},\,u_{i}^{\r}\right)&=&
\left(
  \begin{array}{cc}
u_{i,j}^{\r} & n_{i}^{\r}e_j^{\top} \\
h'_i\left(n_{i}^{\r}\right)e_j&u_{i,j}^{\r}\mbox{\bf I}_d \\
  \end{array}
\right),\,
\end{eqnarray*}
where $(e_1,\cdots,\,e_d)$ is the canonical basis of $\mathbb{R}^d$ and $\mbox{\bf I}_d$
is the $d\times d$ unit matrix, thus \eqref{E-M3} can be written as
\begin{equation}\label{jh}
 \p W_{\v}^{\r}+\sum_{j=1}^d A_\v^j\left(n_{\v}^{\r},\,u_{\v}^{\r}\right)\partial_{x_j}W_{\v}^{\r}
=-H_{\v,\r}^{1} +H_{\v,\r}^{2}-R_{\v}^{\r},\quad \v=e,\,i.
 \end{equation}
with the  initial data
\begin{equation}
\label{jhcz1}
t=0\,:\quad W_{ \v}^{\r}=W_{ \v,0}^{\r},\quad \v=e,\,i,
\end{equation}
in which
\begin{eqnarray*}
   W_{ e,0}^{\r}&=& \left(\begin{matrix}
     N_{e}^{\r}(0,\,\cdot) \\
     {\r}U_{e}^{\r}(0,\,\cdot) \\
   \end{matrix}\right)
= \left(\begin{matrix}
     n_{e,0}^\r-n_{e,\r}^m(0,\,\cdot) \\
      {\r}\left(u_{e,0}^\r-u_{e,\r}^m(0,\,\cdot)\right) \\
   \end{matrix}\right), \nonumber\\
   W_{ i,0}^{\r}&=& \left(\begin{matrix}
     N_{i}^{\r}(0,\,\cdot) \\
     U_{i}^{\r}(0,\,\cdot) \\
   \end{matrix}\right)
= \left(\begin{matrix}
     n_{i,0}^\r-n_{i,\r}^m(0,\,\cdot) \\
     u_{i,0}^\r-u_{i,\r}^m(0,\,\cdot) \\
   \end{matrix}\right).
\end{eqnarray*}
System \eqref{jh}-\eqref{jhcz1} for $W_{ \v}^{\r}$ is symmetrizable hyperbolic when $n_{\v}^{\r}>0$\replaced{, which is ensured by Proposition 2.1.}{. Indeed, since the density $n^0$ of the leading profile satisfies}
\replaced{Indeed}{With this}, let
$$   A_e^0\left(n_{e}^{\r}\right)=
\left(\begin{array}{cc}
h'_e\left(n_{e}^{\r}\right) & 0 \\
0& n_{e}^{\r}\mbox{\bf I}_d \\
  \end{array}\right),
  A_i^0\left(n_{i}^{\r}\right)=
\left(\begin{array}{cc}
h'_i\left(n_{i}^{\r}\right) & 0 \\
0&n_{i}^{\r}\mbox{\bf I}_d \\
  \end{array}\right),
$$
\replaced{thus}{and} for \replaced{$1\leq j\leq d$}{$j=1,2,3$},
\begin{eqnarray*}
\tilde{A}_e^j\left(n_{e}^{\r},\,u_{e}^{\r}\right)
&=&A_e^0\left(n_{e}^{\r}\right)A_e^j\left(n_{e}^{\r},\,u_{e}^{\r}\right)
= \left(
  \begin{array}{cc}
h'_e\left(n_{e}^{\r}\right)u_{e,j}^{\r} &\frac{1}{\r}p'_e\left(n_{e}^{\r}\right) e_j^{\top} \\
\frac{1}{\r}p'_e\left(n_{e}^{\r}\right)e_j&  n_{e}^{\r}u_{e,j}^{\r}\mbox{\bf I}_d \\
  \end{array}
\right),    \nonumber\\
\tilde{A}_i^j\left(n_{i}^{\r},\,u_{i}^{\r}\right)
&=&A_i^0\left(n_{i}^{\r}\right)A_i^j\left(n_{i}^{\r},\,u_{i}^{\r}\right)
=\left(
  \begin{array}{cc}
h'_i\left(n_{i}^{\r}\right)u_{i,j}^{\r} &p'_i\left(n_{i}^{\r}\right) e_j^{\top} \\
p'_i\left(n_{i}^{\r}\right)e_j&n_{i}^{\r}u_{i,j}^{\r}\mbox{\bf I}_d \\
  \end{array}
\right),
\end{eqnarray*}
then $A_\v^0$ is positively definite and $\tilde{A}_\v^j$ is
symmetric for all $1 \leq j \leq d$. Thus, the \replaced{theories}{theorem} of Kato for the local
existence of smooth solutions can also be applied to \eqref{jh}-\eqref{jhcz1}.

By standard arguments, to prove Theorem \ref{e2.1}, it suffices to establish uniform estimates of $W_{\added{\nu}}^{\r}$ with respect to $\r$. \added{Since \eqref{jh}-\eqref{jhcz1} is symmetrizable hyperbolic, we deduce that there exists a time $\bar{T}^\eps\leq T_2^{\eps,1}$ (For convenience, we still denote $\bar{T}^\eps$ as $T_2^{\eps,1}$), such that} \deleted{Since} $W_{\added{\nu}}^{\r}\in C\left(\left[0,\,T_2^{\r,1}\right];\,H^s\left(\mathbb{R}^d\right)\right)$, \added{which implies that} the function $t\rightarrow \|W_{\added{\nu}}^{\r}\|_s$ is continuous on  $\left[0,\,T_2^{\r,1}\right]$. From \eqref{e2.12} and $m \geq 1$, there exists $\tilde{T}^{\r,1}\in\left(0,\,T_2^{\r,1}\right]$ such that
\begin{equation*}
\label{3.6}
    \|W_{\added{\nu}}^{\r}\added{(t)}\|_s \leq C,\quad\forall\,t\in\left[0,\,\tilde{T}^{\r,1}\right],
\end{equation*}
provided that $\r> 0$ is bounded by a constant. If \replaced{$s> \dfrac{d}{2}+1$}{$s \geq 3$}, the embedding from $H^s\left(\mathbb{R}^d\right)$ to $W^{1,\,\infty}\left(\mathbb{R}^d\right)$ is continuous. Then we have
\begin{equation*}
  \|W_{\added{\nu}}^{\r}\added{(t)}\|_{W^{1,\,\infty}\left(\mathbb{R}^d\right)}\leq C,
\quad\forall\,t\in\left[0,\,\tilde{T}^{\r,1}\right].
\end{equation*}
\added{Since $\tilde{T}^{\eps,1}\leq T_2^{\eps,1}\leq T_1^e$, hence for $\eps<1$, it is obvious that the approximate solution satisfies}
\begin{equation*}\label{uniformboundedness}
\added{ \|(n_{\nu,\eps}^m, u_{\nu,\eps}^m,\phi_\eps^m )(t)\|_{W^{1,\,\infty}\left(\mathbb{R}^d\right)}\leq C,
\quad\forall\,t\in\left[0,\,\tilde{T}^{\r,1}\right],}
\end{equation*}
\added{which implies that
\[
 \|(n_e^\eps, n_i^\eps, \eps u_e^\eps, u_i^\eps, \phi^\eps)(t)\|_{W^{1,\,\infty}\left(\mathbb{R}^d\right)}\leq C,
\quad\forall\,t\in\left[0,\,\tilde{T}^{\r,1}\right].
\]
Now it remains to show that $u_e^\eps$ is uniformly bounded with respect to $\eps$. Noticing the initial condition \eqref{e2.12}, since $2m> s$, we immediately have for $\eps<1$,
\[
\|u_{e,0}^{\eps,1}\|_s\leq \|u_{e,0}^{\eps,1}-u_{e,\eps}^{e,m}(0, \cdot)\|_s +\|u_{e,\eps}^{e,m}(0, \cdot)\|_s\leq C_1 \eps^{2m}+C\leq C,
\]
which, by Proposition 2.1, implies that there exists a time $\tilde{T}$, such that $u_e^\eps$ is uniformly bounded with respect to $\eps$ on the time inteval $[0,\tilde{T}]$. This implies that $U_e^\eps$ is also uniformly bounded with respect to $\eps$ on the same time interval. In the following, we will establish the estimates for $W_\nu^\eps$ on the time interval $[0, T^{\r,1}]$, with $T^{\r,1}=\min(\tilde{T}^{\r,1}, \tilde{T})$.}
In order to prove $T_1^{\r,1}\geq T_1^e$, we need to show that there exists a constant
$\mu > 0$ such that
\begin{equation*}
\sup_{0\leq t\leq T^{\r,1}}\|W_{\added{\nu}}^{\r}\added{(t)}\|_s \leq C\r^\mu.
\end{equation*}

\subsubsection{\bf $L^2$-estimates.}
In what follows, we always assume that the conditions of Theorem \ref{e2.1} hold.

\begin{lemma}
\label{lem1}
For all $t\in\left[0,\,T^{\r,1}\right]$ and sufficiently small $\r>0$, we have
\begin{eqnarray}\label{F3.78}
&\qquad&\frac{d}{dt}\left(\sum\limits_{\v=e,\,i}\big\langle A_{\v}^0\left(n_{\v}^{\r}\right)W_{\v}^{\r},\,W_{\v}^{\r}
\big\rangle+\|\nabla\Phi^\r\|^2\right) \nonumber\\
&\leq &C\sum\limits_{\v=e,\,i}\left\|W_{\v}^{\r}\right\|^2+\added{C}\|\nabla\Phi^\r\|^2+C{\r}^{4m\added{+2}}.
\end{eqnarray}
\end{lemma}

\noindent\textbf{Proof.} \underline{Step}1:
Taking the inner product of the \replaced{equations for electrons}{electron equations} in \eqref{jh} with $2A_e^0\left(n_{e}^{\r}\right)W_{e}^{\r}$ in \replaced{$L^2(\mathbb{R}^d)$}{$L^2(\mathbb{R}^3)$}, we obtain the following energy equality for $W_{e}^{\r}$
\begin{align}
\label{F3.8}
    \frac{d}{dt}\big\langle A_e^0(n^\r_{e})W_{e}^{\r},W_{e}^{\r}\big\rangle
=& \big\langle {\mathrm{div}} A_e(n^\r_{e},u^\r_{e})W_{e}^{\r},W_{e}^{\r}\big\rangle
-2\left<A_e^0\left(n_{e}^{\r}\right)W_{e}^{\r},\,H_{e,\r}^{1}\right>
\nonumber\\
&+2\left<A_e^0\left(n_{e}^{\r}\right)W_{e}^{\r},\,H_{e,\r}^{2}\right>
-2\left<A_e^0\left(n_{e}^{\r}\right)W_{e}^{\r},\,R_{e}^{\r}\right>,
\end{align}
where
\begin{equation}\label{divA_e}
\div A_e\left(n_{e}^{\r},\,u_{e}^{\r}\right)=\p A_e^0\left(n_{e}^{\r}\right)+\sum_{j=1}^d \partial_{x_j} \tilde{A}_e^j\left(n_{e}^{\r},\,u_{e}^{\r}\right).
\end{equation}
Now we deal with the right hand side of (\ref{F3.8}) term by term. From the mass conservation equation
 $ \pt n^\r_{e} = - {\mathrm{div}}(n^\r_{e} u^\r_{e}) $, obviously we have
\added{\[
	\|\pt A_e^0(n_e^\eps)\|_\infty\leq C\|\pt n^\r_{e}\|_\infty\leq C\|{\mathrm{div}}(n^\r_{e} u^\r_{e})\|_{s-1}\leq C, 
	\]which implies}\begin{eqnarray}
\label{F3.9}
  \big\langle\pt A_e^0(n^\r_{e})W_{e}^{\r},W_{e}^{\r}\big\rangle
&\leq& C\|W_{e}^{\r}\|^2.
\end{eqnarray}
\replaced{In}{then in} view of the expression of  $ \tilde{A}_i^j(W_{e}^{\r}) $,  \added{for $1\leq j\leq d$,} we obtain
\begin{eqnarray*}
  \big<\pa_{x_j}\tilde{A}_e^j(n_{e}^{\r},\,u_{e}^{\r})W_{e}^{\r},W_{e}^{\r}\big>
&\!\!=\!\!& \big<\pa_{x_j}(h'_{e}(n^\r_{e})u^\r_{e,j})N_{e}^{\r},N_{e}^{\r}\big>
+2\big<N_{e}^{\r}\pa_{x_j}(p'_{e}(n^\r_{e})e_j),{U_{e}^{\r}}\big>\nonumber\\[2mm]
&\!\!\,\!\!& +\, \r^2\big<\pa_{x_j}(n^\r_{e} u^\r_{e,j}){U_{e}^{\r}},{U_{e}^{\r}}\big>,
\end{eqnarray*}
in which we have
\[  \big<\pa_{x_j}(h'_{e}(n^\r_{e})u^\r_{e,j})N_{e}^{\r},N_{e}^{\r}\big>
+ \r^2\big<\pa_{x_j}(n^\r_{e} u^\r_{e,j}){U_{e}^{\r}},{U_{e}^{\r}}\big>
\leq C\|W_{e}^{\r}\|^2,  \]
and
\[ 2\sum_{j=1}^d\big<N_{e}^{\r}\pa_{x_j}(p'_{e}(n^\r_{e})e_j),{U_{e}^{\r}}\big>
= 2 \big<N_{e}^{\r}\D p'_{e}(n^\r_{e}),{U_{e}^{\r}}\big>,  \]
therefore,
\begin{equation}
\label{F3.10}
 \sum_{j=1}^d\big\langle\pa_{x_j}\tilde{A}_e^j(n_{e}^{\r},\,u_{e}^{\r})W_{e}^{\r},W_{e}^{\r}\big\rangle
\leq C\|W_{e}^{\r}\|^2
+   \added{2}\big<N_{e}^{\r}\D p'_{e}(n^\r_{e}),{U_{e}^{\r}}\big>.
\end{equation}
It follows from   \eqref{F3.9} and \eqref{F3.10} that
\begin{equation}
\label{F3.11}
   \big\langle {\mathrm{div}} A_e(n^\r_{e},u^\r_{e})W_{e}^{\r},W_{e}^{\r}\big\rangle
\leq C\|W_{e}^{\r}\|^2
+ 2 \big<N_{e}^{\r}\D p'(n^\r_{e}),{U_{e}^{\r}}\big>.
\end{equation}

For the remaining terms without $H_{e,\eps}^2$ \replaced{on}{in} the right hand side of (\ref{F3.8}), it can be treated as
\begin{eqnarray*}
&\,&-2\big\langle A_e^0(n^\r_{e})W_{e}^{\r},H_{e,\r}^{1}\big\rangle
-2\big\langle A_e^0\left(n_{e}^{\r}\right)W_{e}^{\r},\,R_{e}^{\r}\big\rangle
\nonumber\\
&=& - 2\big\langle h_e'(n^\r_{e})
\big(N_{e}^{\r}{\mathrm{div}} u_{e,\r}^m + {U_{e}^{\r}}\D n_{e,\r}^m
+R_{n_e}^{\r,m}\big),N_{e}^{\r}\big\rangle
 \nonumber\\
&&-\, 2\Big\langle n^\r_{e}\Big[\r^2\big({U_{e}^{\r}}\cdot\D\big) u_{e,\r}^m+\left(h'_{e}(N_e^\r+n_{e,\r}^{m})-h'_{e}(n_{e,\r}^{m})\right)\D n_{e,\r}^{m}
+ R_{u_e}^{\r,m}\Big],{U_{e}^{\r}}\Big\rangle\nonumber\\
&=& -\,2\big\langle h_e'(n^\r_{e})N_{e}^{\r}{\mathrm{div}} u_{e,\r}^m,N_{e}^{\r}\big\rangle
- 2\r^2\big\langle n^\r_{e}\big({U_{e}^{\r}}\cdot\D\big) u_{e,\r}^m,{U_{e}^{\r}}\big\rangle
- 2\big\langle h_e'(n^\r_{e})R_{n_e}^{\r,m},N_{e}^{\r}\big\rangle\nonumber\\
&\,&- 2\big\langle n^\r_{e} R_{u_e}^{\r,m},{U_{e}^{\r}}\big\rangle
-\,2\big\langle h_e'(n^\r_{e})N_{e}^{\r}\D n_{e,\r}^m,{U_{e}^{\r}}\big\rangle
- 2 \big\langle n^\r_{e}\big(h_e'(n^\r_{e})-h_e'(n_{e,\r}^m)\big)
\D n_{e,\r}^m,{U_{e}^{\r}}\big\rangle \\
&\leq & -\,2\big\langle h_e'(n^\r_{e})N_{e}^{\r}\D n_{e,\r}^m,{U_{e}^{\r}}\big\rangle
- 2 \big\langle n^\r_{e}\big(h_e'(n^\r_{e})-h_e'(n_{e,\r}^m)\big)
\D n_{e,\r}^m,{U_{e}^{\r}}\big\rangle \\
&\,&+C\|W_{e}^{\r}\|^2+ C\Big( \|R_{n_e}^{\r,m}\|^2 + \added{\dfrac{1}{\r^2}}\|R_{u_e}^{\r,m}\|^2\Big)\\
&\added{\leq}&\added{ -\,2\big\langle h_e'(n^\r_{e})N_{e}^{\r}\D n_{e,\r}^m,{U_{e}^{\r}}\big\rangle
- 2 \big\langle n^\r_{e}\big(h_e'(n^\r_{e})-h_e'(n_{e,\r}^m)\big)
\D n_{e,\r}^m,{U_{e}^{\r}}\big\rangle} \\
&\,& \added{+C\|W_{e}^{\r}\|^2+C\eps^{4m\added{+2}}.}
\end{eqnarray*}
For the term containing $H_{e,\r}^{2}$ \replaced{on}{in} the right hand side of \eqref{F3.8}, a direct calculation gives
\begin{align*}\label{F3.15}
2\big\langle A_e^0\left(n_{e}^{\r}\right)W_{e}^{\r},\,H_{e,\r}^{2}\big\rangle
=2\big\langle n_{e}^{\r}U_{e}^{\r},\,\nabla\Phi^\r\big\rangle.
\end{align*}
Back to \eqref{F3.8}, combining the above three estimates yield
\begin{eqnarray}
\label{F3.16}
&\qquad& \frac{d}{dt} \big\langle A_e^0(n^\r_{e})W_{e}^{\r},W_{e}^{\r}\big\rangle \nonumber\\
  &\leq& C\|W_{e}^{\r}\|^2+2\big\langle n_{e}^{\r}U_{e}^{\r},\,\nabla\Phi^\r\big\rangle+ C\eps^{4m\added{+2}}+ 2r^\r,
\end{eqnarray}
where the remaining term
\[  r^\r=\big\langle N_{e}^{\r}\D p_e'(n^\r_{e})-h_e'(n^\r_{e})N_{e}^{\r}\D n_{e,\r}^m
+n^\r_{e}\big(h_e'(n_{e,\r}^m)-h_e'(n^\r_{e})\big)\D n_{e,\r}^m,{U_{e}^{\r}}\big\rangle. \]
\replaced{Noticing}{Note} that
\begin{equation}
\label{F3.17}
   n^\r_{e}= N_{e}^{\r}+ n_{e,\r}^m, \quad p_e''(n_e^{\added{\eps}})= h_e'(n_e^{\added{\eps}}) + n_e^{\added{\eps}}h''(n_e^{\added{\eps}}),
\end{equation}
\replaced{and by the Taylor's formula}{When  $ N_{e}^{\r} $  is small}, we have
\begin{equation}
\label{F3.18}
 h_e'(n_{e,\r}^m)-h_e'(n^\r_{e})=-h_e''(n^\r_{e})N_{e}^{\r}+\frac{1}{2}h_e'''(n^\r_{e}-\theta N_{e}^{\r})(N_{e}^{\r})^2,
\mbox{ with }\theta\in[0,1],
\end{equation}
then,
\begin{eqnarray*}
&\,&N_{e}^{\r}\D p_e'(n^\r_{e})-h_e'(n^\r_{e})N_{e}^{\r}\D n_{e,\r}^m +n^\r_{e}\big(h_e'(n_{e,\r}^m)
-h_e'(n^\r_{e})\big)\D n_{e,\r}^m \nonumber\\
&=& N_{e}^{\r} p_e''(n^\r_{e}) \D N_{e}^{\r} + N_{e}^{\r} p_e''(n^\r_{e}) \D n_{e,\r}^m
-h_e'(n^\r_{e})N_{e}^{\r}\D n_{e,\r}^m-n^\r_{e} h_e''(n^\r_{e})N_{e}^{\r}\D n_{e,\r}^m\nonumber\\
&\,&+\frac{1}{2}n^\r_{e} h_e'''(n^\r_{e}-\theta N_{e}^{\r})(N_{e}^{\r})^2\D n_{e,\r}^m \nonumber\\
&=& \frac{1}{2}\,p_e''(n^\r_{e}) \D(N_{e}^{\r})^2+\frac{1}{2}n^\r_{e} h_e'''(n^\r_{e}-\theta N_{e}^{\r})
(N_{e}^{\r})^2\D n_{e,\r}^m.
\end{eqnarray*}
Therefore,
\begin{eqnarray*}
 r^\r&=&\frac{1}{2}\big\langle p_e''(n^\r_{e})\D(N_{e}^{\r})^2,{U_{e}^{\r}}\big\rangle
+\frac{1}{2}\big\langle n^\r_{e} h_e'''(n^\r_{e}-\theta N_{e}^{\r})
(N_{e}^{\r})^2\D n_{e,\r}^m,{U_{e}^{\r}}\big\rangle \nonumber\\
&=& - \frac{1}{2}\big\langle (N_{e}^{\r})^2,{\mathrm{div}}\big(p_e''(n^\r_{e}){U_{e}^{\r}}\big)\big\rangle
+\frac{1}{2}\big\langle n^\r_{e} h_e'''(n^\r_{e}-\theta N_{e}^{\r})
(N_{e}^{\r})^2\D n_{e,\r}^m,{U_{e}^{\r}}\big\rangle \nonumber\\
&\leq& C\|N_{e}^{\r}\|^2.
\end{eqnarray*}
This together  with (\ref{F3.16}) yields
\begin{equation}
\label{F3.19}
\frac{d}{dt} \big\langle A_e^0(n^\r_{e})W_{e}^{\r},W_{e}^{\r}\big\rangle
\leq C\|W_{e}^{\r}\|^2+2\big\langle n_{e}^{\r}U_{e}^{\r},\,\nabla\Phi^\r\big\rangle
+ C\eps^{4m\added{+2}}.
\end{equation}

\underline{Step}2:
Taking the inner product of the \replaced{equations for ions}{ion equations} in \eqref{jh} with $2A_i^0\left(n_{i}^{\r}\right)W_{i}^{\r}$ in \replaced{$L^2(\mathbb{R}^d)$}{$L^2(\mathbb{R}^3)$}, we obtain the following energy equality for $W_{i}^{\r}$
\begin{align}
\frac{d}{dt}\big\langle A_i^0\left(n_{i}^{\r}\right)W_{i}^{\r},\,W_{i}^{\r}\big\rangle
=&-2\big\langle A_i^0\left(n_{i}^{\r}\right)W_{i}^{\r},\,H_{i,\r}^{1}\big\rangle
+2\big\langle A_i^0\left(n_{i}^{\r}\right)W_{i}^{\r},\,H_{i,\r}^{2}\big\rangle\nonumber\\
&-2\big\langle A_i^0\left(n_{i}^{\r}\right)W_{i}^{\r},\,R_{i}^{\r}\big\rangle
+\big\langle\div A_i\left(n_{i}^{\r},\,u_{i}^{\r}\right)W_{i}^{\r},\,W_{i}^{\r}\big\rangle,
\label{16}
\end{align}
where
\begin{equation}\label{divA_i}\div A_i\left(n_{i}^{\r},\,u_{i}^{\r}\right)=\p A_i^0\left(n_{i}^{\r}\right)+\sum_{j=1}^d \partial_{x_j} \tilde{A}_i^j\left(n_{i}^{\r},\,u_{i}^{\r}\right),\end{equation}
which are treated term by term as follows. \replaced{Noticing}{Notice} the expressions of $A_i^0,\,\div A_i$ and $H_{i,\r}^{1}$, we have
\begin{eqnarray*}
\label{17}
\left|\big\langle A_i^0\left(n_{i}^{\r}\right)W_{i}^{\r},\,H_{i,\r}^{1}\big\rangle\right|&\leq& C\left\|W_{i}^{\r}\right\|^2,\\
\label{18}
\left|\big\langle A_i^0\left(n_{i}^{\r}\right)W_{i}^{\r},\,R_{i}^{\r}\big\rangle\right|&\leq& C\left\|W_{i}^{\r}\right\|^2+C\left\|R_{i}^{\r}\right\|^2\leq C\left\|W_{i}^{\r}\right\|^2+C{\r}^{4m\added{+4}},\\
\label{19}
\left|\big\langle\div A_i\left(n_{i}^{\r},\,u_{i}^{\r}\right) W_{i}^{\r},\,W_{i}^{\r}\big\rangle\right|&\leq &C\big\|\div A_i\left(n_{i}^{\r},\,u_{i}^{\r}\right)\big\|_\infty\left\|W_{i}^{\r}\right\|^2\leq C\left\|W_{i}^{\r}\right\|^2.
\end{eqnarray*}
For the term containing $H_{i,\r}^{2}$ in the right hand side of \eqref{16}, a direct calculation gives
\begin{align*}\label{20}
2\big\langle A_i^0\left(n_{i}^{\r}\right)W_{i}^{\r},\,H_{i,\r}^{2}\big\rangle
=-2\big\langle n_{i}^{\r}U_{i}^{\r},\,\nabla\Phi^\r\big\rangle.
\end{align*}
Inserting the above four estimates into \eqref{16}, we get
\begin{equation}
\label{lm1}
\frac{d}{dt}\big\langle A_i^0\left(n_{i}^{\r}\right)W_{i}^{\r},\,W_{i}^{\r}\big\rangle
\leq-2\big\langle n_{i}^{\r}U_{i}^{\r},\,\nabla\Phi^\r\big\rangle+C\left\|W_{i}^{\r}\right\|^2+C{\r}^{4m\added{+4}} .
\end{equation}

\underline{Step}3:
Summing \eqref{F3.19} and \eqref{lm1} \deleted{for all $|\a|\leq s$}, we obtain
\begin{eqnarray*}
\label{l1m1}
\frac{d}{dt}\sum\limits_{\v=e,\,i}\left< A_i^0\left(n_{\v}^{\r}\right)
W_{\v}^{\r},\,W_{\v}^{\r}
\right>&\leq&2\big\langle n_{e}^{\r}U_{e}^{\r},\,\nabla\Phi^\r\big\rangle
-2\big\langle n_{i}^{\r}U_{i}^{\r},\,\nabla\Phi^\r\big\rangle\nonumber\\
&\qquad&+C\sum\limits_{\v=e,\,i}\left\|W_{\v}^{\r}\right\|^2+C{\r}^{4m\added{+2}}\nonumber\\
&=&2\big\langle \left(n_{e}^{\r}u_{e}^{\r}-n_{e,\r}^{m}u_{e,\r}^{m}\right)
-\left(n_{i}^{\r}u_{i}^{\r}-n_{i,\r}^{m}u_{i,\r}^{m}\right),\,\nabla\Phi^\r\big\rangle\nonumber\\
&\,& +\big\langle N_{e}^{\r}u_{e,\r}^{m}-N_{i}^{\r}u_{i,\r}^{m},\,\nabla\Phi^\r\big\rangle\nonumber\\
&\qquad&+C\sum\limits_{\v=e,\,i}\left\|W_{\v}^{\r}\right\|^2+C{\r}^{4m\added{+2}},
\end{eqnarray*}
in which we have
\begin{eqnarray*}
\label{g1j}
&\,&2\big\langle\left(n_{e}^{\r}u_{e}^{\r}-n_{e,\r}^{m}u_{e,\r}^{m}\right)
-\left(n_{i}^{\r}u_{i}^{\r}-n_{i,\r}^{m}u_{i,\r}^{m}\right),\,\nabla\Phi^\r\big\rangle
 \nonumber\\&=& -2\big\langle\div \left(n_{e}^{\r}u_{e}^{\r}-n_{e,\r}^{m}u_{e,\r}^{m}\right)
-\div \left(n_{i}^{\r}u_{i}^{\r}-n_{i,\r}^{m}u_{i,\r}^{m}\right),\,\Phi^\r\big\rangle
 \nonumber\\&=&2\big\langle\pt\left(N_{e}^{\r}-N_{i}^{\r}\right)-\left(R_{n_e}^{\r,m}-R_{n_i}^{\r,m}\right),\,\Phi^\r\big\rangle
 \nonumber\\&=&-\frac{d}{dt}\|\nabla\Phi^\r\|^2
+2\big\langle\tilde{R}_{n_e}^{\r,m}-\tilde{R}_{n_i}^{\r,m},\,\nabla\Phi^\r\big\rangle
\nonumber\\&\leq&-\frac{d}{dt}\|\nabla\Phi^\r\|^2
+\added{C}\|\nabla\Phi^\r\|^2+C{\r}^{4m\added{+2}},
\end{eqnarray*}
and \added{using Lemma \ref{lem2.3},}
\begin{eqnarray*}
\label{gj}
2 \big\langle N_{e}^{\r}u_{e,\r}^{m}-N_{i}^{\r}u_{i,\r}^{m},\,\nabla\Phi^\r\big\rangle
&=&2\big\langle N_{i}^{\r}\left(u_{e,\r}^{m}-u_{i,\r}^{m}\right),\,\nabla\Phi^\r\big\rangle
+2\big\langle\left(N_{e}^{\r}-N_{i}^{\r}\right)u_{e,\r}^{m},\,\nabla\Phi^\r\big\rangle \nonumber\\
&=&2\big\langle N_{i}^{\r}\left(u_{e,\r}^{m}-u_{i,\r}^{m}\right),\,\nabla\Phi^\r\big\rangle
+2\big\langle\triangle\Phi^{\r} u_{e,\r}^{m},\,\nabla\Phi^{\r}\big\rangle
\nonumber\\&\leq& \added{C}\|N_{i}^{\r}\|^2
+\added{C}\|\nabla\Phi^\r\|^2.
\end{eqnarray*}
Combining these estimates yields \eqref{F3.78}. \hfill$\square$

\subsubsection{\bf Higher order estimates.}
Let \replaced{$\alpha\in\N^d$}{$\a\in\mathbb{N}^3$} with $1\leq|\a|\leq s$. Applying $\partial^\a_x$ to \eqref{jh}, we get
\begin{equation}\label{jh2}
\p \partial^\a_xW_{\v}^{\r}+\sum_{j=1}^d A_{\v}^j\left(n_{\v}^{\r},\,u_{\v}^{\r}\right)\partial_{x_j}\partial^\a_xW_{\v}^{\r}
=-\partial^\a_x\left(H_{\v,\r}^{1}-H_{\v,\r}^{2}+R_{\v}^{\r}\right)
+J_{\v,\r}^{\a},\quad \v=e,\,i,
\end{equation}
where
\begin{equation}\label{J}
J_{\v,\r}^{\a}=\sum_{j=1}^d \left(A_{\v}^j\left(n_{\v}^{\r},\,u_{\v}^{\r}\right)\partial_{x_j}\partial^\a_xW_{\v}^{\r}
-\partial^\a_x\left(A_{\v}^j\left(n_{\v}^{\r},\,u_{\v}^{\r}\right)\partial_{x_j}W_{\v}^{\r}\right)\right).
\end{equation}

\begin{lemma}
\label{lem3}
   For all $t\in\left[0,\,T^{\r,1}\right]$ and sufficiently small $\r>0$, we have
\begin{eqnarray}\label{ho}
&\quad&\frac{d}{dt}\left( \sum_{\nu=e,i}\big\langle A_\nu^0(n^\r_{\nu})\pa^\a_xW_{\nu}^{\r},\pa^\a_xW_{\nu}^{\r}\big\rangle+\|\nabla \pa_x^\a\Phi^\eps\|^2\right)\nonumber\\
&\leq& C\sum_{\nu=e,i}\|W_\nu^\eps\|_{|\a|}^2+\added{C}\|\nabla \Phi^\eps\|_{|\a|}^2+\dfrac{C}{\r^2}\|W_e^\r\|_{|\alpha|-1}^2+C\eps^{4m\added{+2}}.
\end{eqnarray}
\end{lemma}
\noindent \textbf{Proof. }\underline{Step}1: Taking the inner product of the \replaced{equations for electrons}{electron equations} in \eqref{jh2} with
$2A_e^0\left(n_{e}^{\r}\right)\partial^\a_xW_{e}^{\r}$ in \replaced{$L^2(\mathbb{R}^d)$}{$L^2(\mathbb{R}^3)$}
yields the following energy equality for $\partial^\a_xW_{e}^{\r}$
\begin{eqnarray}
\label{F3.23}
&\quad&    \frac{d}{dt}\big\langle A_e^0(n^\r_{e})\pa^\a_xW_{e}^{\r},\pa^\a_xW_{e}^{\r}\big\rangle\nonumber\\
&=& \big\langle {\mathrm{div}} A_e(n^\r_{e},u^\r_{e})\pa^\a_xW_{e}^{\r},\pa^\a_xW_{e}^{\r}\big\rangle-2\left<A_e^0\left(n_{e}^{\r}\right)\pa^\a_xW_{e}^{\r},\,\pa^\a_xH_{e,\r}^{1}+\pa^\a_xR_{e}^{\r}\right>\nonumber\\
&\quad&+2\left<A_e^0\left(n_{e}^{\r}\right)\pa^\a_xW_{e}^{\r},\,\pa^\a_xH_{e,\r}^{2}\right>+2\left\langle A_e^0\left(n_{e}^{\r}\right)\partial^\a_xW_{\added{e}}^{\r},\,J_{e,\r}^{\a}\right\rangle,
\end{eqnarray}
\added{in which $\dive A_e$ is defined in \eqref{divA_e}. We will treat the right hand side of the above} \deleted{which are treated} term by term as follows. First, similarly to \eqref{F3.11}, it is easy to get
\begin{equation}
\label{F3.24}
\left|\left\langle{\mathrm{div}} A_e(n^\r_{e},u^\r_{e})\pa^\a_xW_{e}^{\r},\pa^\a_xW_{e}^{\r}\right\rangle\right|
\leq C\left\|W_{e}^{\r}\right\|^2_{|\a|}
+  \added{2}\big<\pa^\a_xN_{e}^{\r}\D p'(n^\r_{e}),\pa^\a_xU_{e}^{\r}\big>.
\end{equation}
For the terms without $H_{e,\eps}^2$ and $J_{e,\eps}^\alpha$ \replaced{on}{in} the right hand side of \eqref{F3.23}, a straightforward calculation yields
\begin{eqnarray*}
\label{F3.25}
&\,& -2\big\langle A_e^0(n^\r_{e})\pa^\a_xW_{e}^{\r},\pa^\a_xH_{e,\r}^{1}+\pa^\a_xR_{e}^{\r}\big\rangle\nonumber\\
&=& - 2\big< h_e'(n^\r_{e}) \big(\pa^\a_x\big(N_{e}^{\r}\dive u_{e,\r}^m \big)
+\pa^\a_xR_{n_{\added{e}}}^{\r,m}\big),\pa^\a_xN_{e}^{\r}\big>- 2\big< h_e'(n^\r_{e})
\pa^\a_xN_{e}^{\r} \D n_{e,\r}^m,\pa^\a_xU_{e}^{\r}\big>\nonumber\\
&\,& - 2\big< h_e'(n^\r_{e})\big(\pa^\a_x\big({U_{e}^{\r}}\D n_{e,\r}^m\big)
-\pa^\a_xU_{e}^{\r} \D n_{e,\r}^m\big),\pa^\a_xN_{e}^{\r}\big>\nonumber\\
&\,& -\, 2\big< n^\r_{e}\pa^\a_x\Big[\r^2\big({U_{e}^{\r}}\cdot\D\big) u_{e,\r}^m
+ R_{u_{\added{e}}}^{\r,m}\Big],\pa^\a_xU_{e}^{\r}\big>
- 2 \left< n^\r_{e}\pa^\a_x\big(h_e'(n^\r_{e})-h_e'(n_{e,\r}^m)\big)
\D n_{e,\r}^m,\pa^\a_xU_{e}^{\r}\right>\nonumber\\
&\,&- 2\big<n^\r_{e} \big[\pa^\a_x\big((h_e'(n^\r_{e})-h_e'(n_{e,\r}^m))
\D n_{e,\r}^m \big)-\pa^\a_x\big(h_e'(n^\r_{e})-h_e'(n_{e,\r}^m)\big)
\D n_{e,\r}^m\big],\pa^\a_xU_{e}^{\r}\big>,
\end{eqnarray*}
in which we have
\begin{eqnarray*}
\label{F3.26}
- 2\left< h_e'(n^\r_{e}) \big(\pa^\a_x(N_{e}^{\r}\dive u_{e,\r}^m )
+\pa^\a_xR_{n_{\added{e}}}^{\r,m}\big),\pa^\a_xN_{e}^{\r}\right>
&\leq & C\Big(\|W_{e}^{\r}\|^2_{|\a|}+\|R_{n_{\added{e}}}^{\r,m}\|^2_{|\a|}\Big),\\
\label{F3.27}
    -\, 2\left< n^\r_{e}\pa^\a_x\Big[\r^2\big({U_{e}^{\r}}\cdot\D\big) u_{e,\r}^m
+ R_{u_{\added{e}}}^{\r,m}\Big],\pa^\a_xU_{e}^{\r}\right>
&\leq &C\Big(\|W_{e}^{\r}\|^2_{|\a|}
+ \added{\frac{1}{\eps^2}}\|R_{u_{\added{e}}}^{\r,m}\|^2_{|\a|}\Big).
\end{eqnarray*}
Besides, applying the Moser-type inequalities, we have
\begin{align*}
  - 2\left< h_e'(n^\r_{e})\big(\pa^\a_x({U_{e}^{\r}}\D n_{e,\r}^m)
-\pa^\a_xU_{e}^{\r} \D n_{e,\r}^m\big),\pa^\a_xN_{e}^{\r}\right>
\leq C\Big(\|N_{e}^{\r}\|^2_{|\a|} + \|{U_{e}^{\r}}\|^2_{|\a|-1}\Big),
\end{align*}
and \added{by the Taylor's formula}
\begin{eqnarray*}
&\quad&-2\left< n^\r_{e} \big(\pa^\a_x(\big(h_e'(n^\r_{e})-h_e'(n_{e,\r}^m)\big)
\D n_{e,\r}^m )-\pa^\a_x\big(h_e'(n^\r_{e})-h_e'(n_{e,\r}^m)\big)
\D n_{e,\r}^m\big),\pa^\a_xU_{e}^{\r}\right>\nonumber\\
&\leq&\quad C\Big(\frac{1}{\r^2}\|N_{e}^{\r}\|^2_{|\a|-1} +\r^2 \|{U_{e}^{\r}}\|^2_{|\a|}\Big).
\end{eqnarray*}
The above four estimates imply
\begin{eqnarray}
&\,& -2\big\langle A_e^0(n^\r_{e})\pa^\a_xW_{e}^{\r},\pa^\a_xH_{e,\r}^{1}+\pa^\a_xR_{e}^{\r}\big\rangle\nonumber\\
&\leq& \frac{C}{\r^2}\Big(\|N_{e}^{\r}\|^2_{|\a|-1} + \r^2\|{U_{e}^{\r}}\|^2_{|\a|-1}\Big)+C\Big(\|R_{n_{\added{e}}}^{\r,m}\|^2_{|\a|} + \added{\frac{1}{\eps^2}}\|R_{u_{\added{e}}}^{\r,m}\|^2_{|\a|}\Big)\nonumber\\
&\quad&+C\Big(\|N_{e}^{\r}\|^2_{|\a|}+ \r^2\|{U_{e}^{\r}}\|^2_{|\a|}\Big)- 2\big< h_e'(n^\r_{e})\pa^\a_xN_{e}^{\r} \D n_{e,\r}^m,\pa^\a_xU_{e}^{\r}\big>\nonumber\\
&\quad&- 2 \left< n^\r_{e}\pa^\a_x\big(h_e'(n^\r_{e})-h_e'(n_{e,\r}^m)\big)
\D n_{e,\r}^m,\pa^\a_xU_{e}^{\r}\right>.
\end{eqnarray}
For the term containing  $ J_{e,\r}^{\a} $  \replaced{on}{in} the right hand side of \eqref{F3.23}, for \replaced{$1\leq j\leq d$}{$j=1,2,3$} \added{and $U_e^\eps=(U_{e,1}^\eps, \cdots, U_{e,d}^\eps)$}, we have
\begin{eqnarray}\label{jam}
&\;& \left<A_e^0(n^\r_{e})\pa^\alpha_xW_e^\r,\pa^\alpha_x (A_e^j(n^\r_{e},u^\r_{e})\pa_{x_j}W_e^\r)
- A_e^j(n^\r_{e},u^\r_{e})\pa^\alpha_x (\pa_{x_j} W_e^\r)\right> \nonumber\\[2mm]
& = & \left<h_e'(n^\r_{e}) \left(\pa^\alpha_x (u^\r_{e,j} \pa_{x_j} N_e^\r)
- u^\r_{e,j} \pa^\alpha_x \pa_{x_j} N_e^\r\right),\pa^\alpha_x N_e^\r\right>
\nonumber\\[2mm]
& \, & + \;  \r^2\left<n^\r_{e} \left(\pa^\alpha_x (u^\r_{e,j} \pa_{x_j} U_e^\eps)
- u^\r_{e,j}\pa^\alpha_x\pa_{x_j} U_e^\eps\right),\pa^\alpha_x U_e^\eps\right> \nonumber\\[2mm]
& \, & + \; \left<n^\r_{e} \left(\pa^\alpha_x(h_e'(n^\r_{e})\pa_{x_j} N_e^\r)
-h_e'(n^\r_{e})\pa^\alpha_x\pa_{x_j} N_e^\r\right),\pa^\alpha_x U_{e,j}^\r\right>\nonumber\\[2mm]
& \, & + \;\left<h_e'(n^\r_{e})\left(\pa^\alpha_x(n^\r_{e} \pa_{x_j} U_{e,j}^\r)
- n^\r_{e} \pa^\alpha_x \pa_{x_j} U_{e,j}^\r \right),\pa^\alpha_x N_e^\r\right>\nonumber\\[2mm]
& = & \left<h_e'(n^\r_{e}) \left(\pa^\alpha_x (u^\r_{e,j} \pa_{x_j} N_e^\r)
- u^\r_{e,j} \pa^\alpha_x \pa_{x_j} N_e^\r\right),\pa^\alpha_x N_e^\r\right>
\nonumber\\[2mm]
& \, & + \;  \r^2\left<n^\r_{e} \left(\pa^\alpha_x (u^\r_{e,j} \pa_{x_j} U_e^\eps)
- u^\r_{e,j}\pa^\alpha_x\pa_{x_j} U_e^\eps\right),\pa^\alpha_x U_e^\eps\right> \nonumber\\[2mm]
& \, & + \; \left<n^\r_{e} \left(\pa^\alpha_x(h_e'(n^\r_{e})\pa_{x_j} N_e^\r)
-h_e'(n^\r_{e})\pa^\alpha_x\pa_{x_j} N_e^\r-\sum_{{\substack{1\leq i\leq d\\\alpha_i\neq 0}}  }\alpha_i\pa_{x_i} h_e'(n^\r_{e})\pa^{\alpha^i}_x \pa_{x_j} N_e^\r\right),\pa^\alpha_x U_{e,j}^\r\right>\nonumber\\[2mm]
& \, & + \;\left<h_e'(n^\r_{e})\left(\pa^\alpha_x(n^\r_{e} \pa_{x_j} U_{e,j}^\r)
- n^\r_{e} \pa^\alpha_x \pa_{x_j} U_{e,j}^\r -\sum_{{\substack{1\leq i\leq d\\\alpha_i\neq 0}}  }\alpha_i\pa_{x_i}n^\r_{e}\pa^{\alpha^i}_x \pa_{x_j} U_{e,j}^\r\right),\pa^\alpha_x N_e^\r\right>\nonumber\\[2mm]
& \, & + \; \sum_{{\substack{1\leq i\leq d\\\alpha_i\neq 0}}  }\alpha_i
\left(\left<n^\r_{e}\pa_{x_i}h_e'(n^\r_{e})\pa^{\alpha^i}_x \pa_{x_j} N_e^\r
,\pa^\alpha_x U_{e,j}^\r\right>+
\left<h_e'(n^\r_{e})\pa_{x_i}n^\r_{e}\pa^{\alpha^i}_x \pa_{x_j} U_{e,j}^\r
,\pa^\alpha_x N_e^\r\right>\right)\nonumber\\[2mm]
& \leq & C\left(\|N_e^\r\|^2_{|\a|} + \r^2 \|U_e^\r\|^2_{|\a|}\right)+\frac{C}{\r^2}\left(\|N_e^\r\|^2_{|\a|-1} + \r^2\|U_e^\r\|^2_{|\a|-1}\right)\nonumber\\[2mm]
& \, & + \; \sum_{{\substack{1\leq i\leq d\\\alpha_i\neq 0}}  }\alpha_i\left(\left<n^\r_{e}\pa_{x_i}
h_e'(n^\r_{e})\pa^{\alpha^i}_x \pa_{x_j} N_e^\r
,\pa^\alpha_x U_{e,j}^\r\right>+\left<h_e'(n^\r_{e})\pa_{x_i}n^\r_{e}\pa^{\alpha^i}_x \pa_{x_j} U_{e,j}^\r
,\pa^\alpha_x N_e^\r\right>\right),\nonumber
\end{eqnarray}
where  $ \alpha^i $  is also a multi-index and  $\pa_{x_i}\pa^{\alpha^i}_x=\pa^{\alpha}_x  $. The last \added{two} term{\added{s}} \deleted{of \eqref{jam}}  can be estimated as \added{follows. Taking integration by parts with respective to $x_i$ and $x_j$ successively to the first of the last two terms, we obtain}
\begin{eqnarray*}
&\;& \sum_{{\substack{1\leq i\leq d\\\alpha_i\neq 0}}  }\alpha_i\left(\left<n^\r_{e}\pa_{x_i}h_{\added{e}}'(n^\r_{e})\pa^{\alpha^i}_x \pa_{x_j} N_e^\r
,\pa^\alpha_x U_{e,j}^\r\right>+\left<h_{\added{e}}'(n^\r_{e})\pa_{x_i}n^\r_{e}\pa^{\alpha^i}_x \pa_{x_j} U_{e,j}^\r
,\pa^\alpha_x N_e^\r\right>\right) \nonumber\\[2mm]
& = &\sum_{{\substack{1\leq i\leq d\\\alpha_i\neq 0}}  }\alpha_i\left(-\left<n^\r_{e}\pa_{x_i}h_{\added{e}}'(n^\r_{e})\pa^{\alpha}_x \pa_{x_j} N_e^\r
,\pa^{\alpha^i}_x U_{e,j}^\r\right>+\left<h_{\added{e}}'(n^\r_{e})\pa_{x_i}n^\r_{e}\pa^{\alpha^i}_x \pa_{x_j} U_{e,j}^\r
,\pa^\alpha_x N_e^\r\right>\right)\nonumber\\[2mm]
& \, & - \;\sum_{{\substack{1\leq i\leq d\\\alpha_i\neq 0}}  }\alpha_i\left<\pa_{x_i}\left(n^\r_{e}\pa_{x_i}h_{\added{e}}'(n^\r_{e})\right)\pa^{\alpha^i}_x \pa_{x_j} N_e^\r,\pa^{\alpha^i}_x U_{e,j}^\r\right>\nonumber\\[2mm]
& = &
\sum_{{\substack{1\leq i\leq d\\\alpha_i\neq 0}}  }\alpha_i\left(\left<n^\r_{e}\pa_{x_i}h_{\added{e}}'(n^\r_{e})\pa^{\alpha}_x N_e^\r
,\pa^{\alpha^i}_x \pa_{x_j} U_{e,j}^\r\right>+\left<h_{\added{e}}'(n^\r_{e})\pa_{x_i}n^\r_{e}\pa^{\alpha^i}_x \pa_{x_j} U_{e,j}^\r
,\pa^\alpha_x N_e^\r\right>\right)\nonumber\\[2mm]
& \, & - \;\sum_{{\substack{1\leq i\leq d\\\alpha_i\neq 0}}  }\alpha_i\left(\left<\pa_{x_i}\left(n^\r_{e}\pa_{x_i}h_{\added{e}}'(n^\r_{e})\right)\pa^{\alpha^i}_x \pa_{x_j} N_e^\r
,\pa^{\alpha^i}_x U_{e,j}^\r\right>-\left<\pa_{x_j}\left(n^\r_{e}\pa_{x_i}h_{\added{e}}'(n^\r_{e})\right)\pa^{\alpha}_x N_e^\r
,\pa^{\alpha^i}_x U_{e,j}^\r\right>\right)\nonumber\\[2mm]
& \leq & \sum_{{\substack{1\leq i\leq d\\\alpha_i\neq 0}}  }\alpha_i\left<\pa_{x_i}p_{\added{e}}'(n^\r_{e})\pa^{\alpha^i}_x \pa_{x_j} U_{e,j}^\r
,\pa^\alpha_x N_e^\r\right>+ C\left(\|N_e^\r\|^2_{|\a|} + \|U_e^\r\|^2_{|\a|-1}\right),
\end{eqnarray*}
\replaced{where we have used the definition of the enthalpy function. Noting that the fluid is isothermal, i.e.,}{Note that}  $ p_e'(n^\r_{e})=a_e^2 $,  we get
\begin{equation}
\label{e220}
\left|\left\langle A_e^0\left(n^{\r}\right)\partial^\a_xW_e^\r,\,J_{e,\r}^{\a}\right\rangle\right|
\leq C \|W_e^\r\|^2_{|\a|}  +\frac{C}{\r^2} \|W_e^\r\|^2_{|\a|-1}  .
\end{equation}
As to the term containing $H_{e,\r}^{2}$ in the right hand side of \eqref{F3.23}, a direct calculation gives
\begin{align*}\label{Fe3.15}
2\big\langle A_e^0\left(n_{e}^{\r}\right)\pa^\a_xW_{e}^{\r},\,\pa^\a_xH_{e,\r}^{2}\big\rangle
=2\big\langle(n_{e}^{\r}\pa^\a_xU_{e}^{\r},\,\nabla\pa^\a_x\Phi^\r\big\rangle.
\end{align*}
Therefore, \eqref{F3.24}-\eqref{e220} yield
\begin{eqnarray}
\label{F3.31}
\frac{d}{dt}\big\langle A_e^0(n^\r_{e})\pa^\a_xW_{e}^{\r},\pa^\a_xW_{e}^{\r}\big\rangle
&\added{\leq}& \added{C \|W_e^\r\|^2_{|\a|}  +\frac{C}{\r^2} \|W_e^\r\|^2_{|\a|-1}+C\eps^{4m\added{+2}}}\nonumber\\
&\quad&\added{ + 2\big\langle n_{e}^{\r}\pa^\a_xU_{e}^{\r},\,\nabla\pa^\a_x\Phi^\r\big\rangle+2r^\r_\a},
\end{eqnarray}
where
\[  r^\r_\a=\left< \pa^\a_xN_{e}^{\r}\D p_e'(n^\r_{e})-h_e'(n^\r_{e})\pa^\a_xN_{e}^{\r}\D n_{e,\r}^m
+n^\r_{e}\pa^\a_x\big(h_e'(n_{e,\r}^m)-h'_e(n^\r_{e})\big)\D n_{e,\r}^m,\pa^\a_xU_{e}^{\r}\right>, \]
which can be estimated in a \replaced{similar}{same} way as \added{that} in the proof of Lemma \ref{lem1}. Indeed,
by \eqref{F3.17}-\eqref{F3.18}, we have
\begin{eqnarray*}
&\,& \pa^\a_xN_{e}^{\r}\D p_e'(n^\r_{e})-h_e'(n^\r_{e})\pa^\a_xN_{e}^{\r}\D n_{e,\r}^m
+n^\r_{e}\pa^\a_x\big(h_e'(n_{e,\r}^m)-h_e'(n^\r_{e})\big)\D n_{e,\r}^m \nonumber\\
&=& \pa^\a_xN_{e}^{\r} p_e''(n^\r_{e}) \D N_{e}^{\r} + \pa^\a_xN_{e}^{\r} p_e''(n^\r_{e}) \D n_{e,\r}^m
-h_e'(n^\r_{e})\pa^\a_xN_{e}^{\r}\D n_{e,\r}^m \nonumber\\
&\,&-n^\r_{e}\pa^\a_x\Big[h_e''(n^\r_{e})N_{e}^{\r}
-\frac{1}{2}h_e'''(n^\r_{e}-\theta N_{e}^{\r})(N_{e}^{\r})^2\Big]\D n_{e,\r}^m
\nonumber\\
&=& \pa^\a_xN_{e}^{\r} p_e''(n^\r_{e}) \D N_{e}^{\r}
-n^\r_{e}\big[\pa^\a_x\big(h_e''(n^\r_{e})N_{e}^{\r}\big)-h_e''(n^\r_{e})\pa^\a_xN_{e}^{\r}\big]\D n_{e,\r}^m \nonumber\\
&\,&+\frac{1}{2}n^\r_{e}\pa^\a_x\big(h_e'''(n^\r_{e}-\theta N_{e}^{\r})(N_{e}^{\r})^2\big)\D n_{e,\r}^m,
\end{eqnarray*}
by the Moser-type inequalities, we have
\begin{eqnarray*}
  &\,& \|\pa^\a_xN_{e}^{\r}\D p_e'(n^\r_{e})-h_e'(n^\r_{e})\pa^\a_xN_{e}^{\r}\D n_{e,\r}^m
+n^\r_{e}\pa^\a_x\big(h_e'(n_{e,\r}^m)-h_e'(n^\r_{e})\big)\D n_{e,\r}^m\| \\
&\leq& C \|N_{e}^{\r}\|_{|\a|}^2 + C\|N_{e}^{\r}\|_{|\a|-1}+C \|N_{e}^{\r}\|_{|\a|}^2 \\
&\leq& C \|N_{e}^{\r}\|_{|\a|}^2 + C\|N_{e}^{\r}\|_{|\a|-1}.
\end{eqnarray*}
It follows that
\begin{eqnarray}
\label{F3.32}
   2r^\r_\a
&\leq& C\big(\|N_{e}^{\r}\|_{|\a|}^2+\|N_{e}^{\r}\|_{|\a|-1}\big)\|{U_{e}^{\r}}\|_{|\a|}\nonumber\\
&\leq &C\Big( \|N_{e}^{\r}\|^2_{|\a|} + \r^2\|{U_{e}^{\r}}\|^2_{|\a|}\Big)
+ \frac{C}{\r^2}\|N_{e}^{\r}\|^2_{|\a|-1}.
\end{eqnarray}
Combining \eqref{F3.31} and \eqref{F3.32}, we obtain
\begin{eqnarray}
\label{F13.41}
&\quad&\frac{d}{dt}\big\langle A_e^0(n^\r_{e})\pa^\a_xW_{e}^{\r},\pa^\a_xW_{e}^{\r}\big\rangle\nonumber\\
&\leq&2\big\langle n_{e}^{\r}\pa^\a_xU_{e}^{\r},\,\nabla\pa^\a_x\Phi^\r\big\rangle+
 \|W_{e}^{\r}\|^2_{|\a|}  +\frac{C}{\r^2} \|W_{e}^{\r}\|^2_{|\a|-1}+C\eps^{4m\added{+2}}.
\end{eqnarray}

\underline{Step}2:
Taking the inner product of the \replaced{equations for ions}{ion equations} in \eqref{jh2} with
$2A_i^0\left(n_{i}^{\r}\right)\partial^\a_xW_{i}^{\r}$ in \replaced{$L^2(\mathbb{R}^d)$}{$L^2(\mathbb{R}^3)$}
yields the following energy equality for $\partial^\a_xW_{i}^{\r}$
\begin{align}
\label{216}
\frac{d}{dt}\big\langle A_i^0\left(n_{i}^{\r}\right)\partial^\a_xW_{i}^{\r},\,\partial^\a_xW_{i}^{\r}\big\rangle
=&-2\big\langle A_i^0\left(n_{i}^{\r}\right)\partial^\a_xW_{i}^{\r},\,\partial^\a_xH_{i,\r}^{1}\big\rangle
+2\big\langle A_i^0\left(n_{i}^{\r}\right)\partial^\a_xW_{i}^{\r},\,\partial^\a_xH_{i,\r}^{2}\big\rangle\nonumber\\
&-2\big\langle A_i^0\left(n_{i}^{\r}\right)\partial^\a_xW_{i}^{\r},\,\partial^\a_xR_{i}^{\r}\big\rangle
+2\big\langle A_i^0\left(n_{i}^{\r}\right)\partial^\a_xW_{i}^{\r},\,J_{i,\r}^{\a}\big\rangle\nonumber\\
&+\big\langle\div A_i\left(n_{i}^{\r},\,u_{i}^{\r}\right)\partial^\a_xW_{i}^{\r},\,\partial^\a_xW_{i}^{\r}\big\rangle,
\end{align}
\added{where $\dive A_i$ is defined in \eqref{divA_i} and $J_{i,\eps}^\a$ is defined in \eqref{J}.} By \eqref{2.11} and \replaced{uniform boundedness of $n_i^\eps$ and $n_e^\eps$}{\eqref{3.6}}, it is clear that
\begin{eqnarray*}
\left|\big\langle A_i^0\left(n_{i}^{\r}\right)\partial^\a_xW_{i}^{\r},\,\partial^\a_xH_{i,\r}^{1}\big\rangle\right|&\leq &C\left\|W_{i}^{\r}\right\|^2_{|\a|},\\
\left|\big\langle A_i^0\left(n_{i}^{\r}\right)\partial^\a_xW_{i}^{\r},\,\partial^\a_xR_{i}^{\r}\big\rangle\right|&\leq & C\left\|W_{i}^{\r}\right\|^2_{|\a|}+C\left\|R_{i}^{\r}\right\|^2_{|\a|}\nonumber\\
&\leq &C\left\|W_{i}^{\r}\right\|^2_{|\a|}+C\eps^{4m\added{+4}},
\end{eqnarray*}
and
\begin{eqnarray*}
\label{219}
\left|\big\langle \div A_i\left(n_{i}^{\r},\,u_{i}^{\r}\right)\partial^\a_xW_{i}^{\r},\,\partial^\a_xW_{i}^{\r}\big\rangle\right|
&\leq& C\big\|\div A_i(n_{i}^{\r},\,u_{i}^{\r})\big\|_\infty
\left\|W_{i}^{\r}\right\|^2_{|\a|}\nonumber\\
&\leq& C\left\|W_{i}^{\r}\right\|^2_{|\a|}.
\end{eqnarray*}
For the term containing $J_{i,\r}^{\a}$ \replaced{on}{in} the right hand side of \eqref{216}, applying \added{the} Moser-type  inequalities to $J_{i,\r}^{\a}$ together with \replaced{uniform boundedness of $n_i^\eps$ and $n_e^\eps$}{\eqref{3.6}} yields
\begin{equation*}
\label{220}
\left|\big\langle A_i^0\left(n_{i}^{\r}\right)\partial^\a_xW_{i}^{\r},\,J_{i,\r}^{\a}\big\rangle\right|\leq C\left\|W_{i}^{\r}\right\|^2_{|\a|}.
\end{equation*}
\deleted{Since concerning the $L^2(\mathbb{R}^3)$ estimate, we may write the term containing
$\partial^\a_xH_{i,\r}^{2}$ as} A direct computation yields
\begin{align*}\label{221}
\big\langle A_i^0\left(n_{i}^{\r}\right)\partial^\a_xW_{i}^{\r},\,\partial^\a_xH_{i,\r}^{2}\big\rangle
=-2\big\langle n_{i}^{\r}\partial^\a_xU_{i}^{\r},\,\nabla\partial^\a_x\Phi^\r\big\rangle.
\end{align*}
Combining the above five estimates and \eqref{216}, we get
\begin{equation}
\label{lm3}
   \frac{d}{dt}\big\langle A_i^0\left(n_{i}^{\r}\right)\partial^\a_xW_{i}^{\r},\,\partial^\a_xW_{i}^{\r}\big\rangle
\leq-2\big\langle n_{i}^{\r}\partial^\a_xU_{i}^{\r},\,\nabla\partial^\a_x\Phi^\r\big\rangle+C\left\|W_{i}^{\r}\right\|^2_{|\a|}+C{\r}^{4m\added{+4}}.
\end{equation}

\underline{Step}3:
Summing \eqref{F13.41} and \eqref{lm3} \deleted{for all $|\a|\leq s$}, we obtain
\begin{eqnarray}
&\,&\frac{d}{dt}\sum\limits_{\v=e,\,i}\big\langle A_{\v}^0\left(n_{\v}^{\r}\right)
\partial^\a_xW_{\v}^{\r},\,\partial^\a_xW_{\v}^{\r}
\big\rangle\nonumber\\
&\leq&2\big\langle n_{e}^{\r}\partial^\a_xU_{e}^{\r},\,\nabla\partial^\a_x\Phi^\r\big\rangle
-2\big\langle n_{i}^{\r}\partial^\a_xU_{i}^{\r},\,\nabla\partial^\a_x\Phi^\r\big\rangle\nonumber\\
&\quad&+C\sum\limits_{\v=e,\,i}\left\|W_{\v}^{\r}\right\|^2_{|\a|}+\frac{C}{\r^2} \|W_{e}^{\r}\|^2_{|\a|-1}+C{\r}^{4m\added{+2}}
\nonumber\\
&=&2\big\langle\partial^\a_x\left(n_{e}^{\r}U_{e}^{\r}\right),\,\nabla\partial^\a_x\Phi^\r\big\rangle
-2\big\langle\partial^\a_x\left(n_{i}^{\r}U_{i}^{\r}\right),\,\nabla\partial^\a_x\Phi^\r\big\rangle\nonumber\\
&\quad&+C\sum\limits_{\v=e,\,i}\left\|W_{\v}^{\r}\right\|^2_{|\a|}+\frac{C}{\r^2} \|W_{e}^{\r}\|^2_{|\a|-1}+C{\r}^{4m\added{+2}}\nonumber\\
&\,&+2\big\langle n_{e}^{\r}\partial^\a_xU_{e}^{\r}-\partial^\a_x\left(n_{e}^{\r}U_{e}^{\r}\right),\,\nabla\partial^\a_x\Phi^\r\big\rangle
-2\big\langle n_{i}^{\r}\partial^\a_xU_{i}^{\r}-\partial^\a_x\left(n_{i}^{\r}U_{i}^{\r}\right),\,\nabla\partial^\a_x\Phi^\r\big\rangle
\nonumber\\
&=&2\big\langle\partial^\a_x\left(n_{e}^{\r}u_{e}^{\r}-n_{e,\r}^{m}u_{e,\r}^{m}\right)
-\partial^\a_x\left(n_{i}^{\r}u_{i}^{\r}-n_{i,\r}^{m}u_{i,\r}^{m}\right),\,\nabla\partial^\a_x\Phi^\r\big\rangle\nonumber\\
&\quad& -2\big\langle\partial^\a_x\left(N_{e}^{\r}u_{e,\r}^{m}-N_{i}^{\r}u_{i,\r}^{m}\right),\,\nabla\partial^\a_x\Phi^\r\big\rangle\nonumber\\
&\,&+C\sum\limits_{\v=e,\,i}\left\|W_{\v}^{\r}\right\|^2_{|\a|}+\added{\frac{C}{\r^2} \|W_{e}^{\r}\|^2_{|\a|-1}}+C{\r}^{4m\added{+2}}
\nonumber\\
&\,&+2\big\langle n_{e}^{\r}\partial^\a_xU_{e}^{\r}-\partial^\a_x\left(n_{e}^{\r}U_{e}^{\r}\right),\,\nabla\partial^\a_x\Phi^\r\big\rangle-2\big\langle n_{i}^{\r}\partial^\a_xU_{i}^{\r}-\partial^\a_x\left(n_{i}^{\r}U_{i}^{\r}\right),\,\nabla\partial^\a_x\Phi^\r\big\rangle, \nonumber
\end{eqnarray}
in which we have
\begin{eqnarray*}
\label{ggg1j}
&\,&2\big\langle\partial^\a_x\left(n_{e}^{\r}u_{e}^{\r}-n_{e,\r}^{m}u_{e,\r}^{m}\right)
-\partial^\a_x\left(n_{i}^{\r}u_{i}^{\r}-n_{i,\r}^{m}u_{i,\r}^{m}\right),\,\nabla\partial^\a_x\Phi^\r\big\rangle
 \nonumber\\&=& -2\big\langle \div \left(\partial^\a_x\left(n_{e}^{\r}u_{e}^{\r}-n_{e,\r}^{m}u_{e,\r}^{m}\right)\right)
-\div\left(\partial^\a_x \left(n_{i}^{\r}u_{i}^{\r}-n_{i,\r}^{m}u_{i,\r}^{m}\right)\right),\,\partial^\a_x\Phi^\r\big\rangle
 \nonumber\\&=&2\big\langle\pt\left(\partial^\a_xN_{e}^{\r}-\partial^\a_xN_{i}^{\r}\right)-\left(\partial^\a_xR_{n_e}^{\r,m}-\partial^\a_xR_{n_i}^{\r,m}\right),\,\partial^\a_x\Phi^\r\big\rangle
 \nonumber\\&=&-\frac{d}{dt}\|\nabla\partial^\a_x\Phi^\r\|^2
+2\big\langle\partial^\a_x\tilde{R}_{n_e}^{\r,m}-\partial^\a_x\tilde{R}_{n_i}^{\r,m},\,\nabla\partial^\a_x\Phi^\r\big\rangle
\nonumber\\&\leq&-\frac{d}{dt}\|\nabla\partial^\a_x\Phi^\r\|^2
+\|\nabla\partial^\a_x\Phi^\r\|^2+C{\r}^{4m\added{+2}},
\end{eqnarray*}
and \added{using Lemma \ref{lem2.3}}, 
\begin{eqnarray*}
\label{gggj}
&\,&2
 \big\langle\partial^\a_x\left(N_{e}^{\r}u_{e,\r}^{m}-N_{i}^{\r}u_{i,\r}^{m}\right),\,\nabla\partial^\a_x\Phi^\r\big\rangle\nonumber\\
&=&2\big\langle \partial^\a_x\left(N_{i}^{\r}\left(u_{e,\r}^{m}-u_{i,\r}^{m}\right)\right),\,\nabla\partial^\a_x\Phi^\r\big\rangle
+2\big\langle\partial^\a_x\left(\left(N_{e}^{\r}-N_{i}^{\r}\right)u_{e,\r}^{m}\right),\,\nabla\partial^\a_x\Phi^\r\big\rangle\nonumber\\
&=&2\big\langle\partial^\a_x\left(N_{i}^{\r}\left(u_{e,\r}^{m}-u_{i,\r}^{m}\right)\right),\,\nabla\partial^\a_x\Phi^\r\big\rangle
+2\big\langle\triangle\partial^\a_x\Phi^{\r} u_{e,\r}^{m},\,\nabla\partial^\a_x\Phi^{\r}\big\rangle\nonumber\\
&\,&
-2\big\langle\triangle\partial^\a_x\Phi^{\r} u_{e,\r}^{m}-\partial^\a_x\left(\triangle\Phi^{\r} u_{e,\r}^{m}\right),\,\nabla\partial^\a_x\Phi^{\r}\big\rangle
\nonumber\\&\leq& \added{C}\|N_{i}^{\r}\|^2_{|\a|}
+\added{C}\|\nabla\Phi^\r\|^2_{|\a|}.
\end{eqnarray*}
Finally, the Moser-type inequalities \deleted{and the fact $n_{\v}=n_{\v,\r}^{m}+N_{\v}^{\r}\leq C,\mbox{ for }\v=i,e$} imply
\begin{eqnarray*}
\label{g352}
&\,&2
 \big\langle n_{e}^{\r}\partial^\a_xU_{e}^{\r}-\partial^\a_x\left(n_{e}^{\r}U_{e}^{\r}\right),\,\nabla\partial^\a_x\Phi^\r\big\rangle
 -2\big\langle n_{i}^{\r}\partial^\a_xU_{i}^{\r}-\partial^\a_x\left(n_{i}^{\r}U_{i}^{\r}\right),\,\nabla\partial^\a_x\Phi^\r\big\rangle
\nonumber\\&\leq& \frac{C}{\r^2} \|W_{e}^{\r}\|^2_{|\a|-1}+C \|W_{i}^{\r}\|^2_{|\a|-1}+\added{C}\|\nabla\Phi^\r\|^2_{|\a|}.
\end{eqnarray*}
Combining the above four inequalities yields \eqref{ho}. \hfill$\square$
%
%

\subsection{\bf Proof of Theorem \ref{e2.1}}
We deal with \eqref{ho} by induction for $ 1\leq |\a|\leq s $. In view of
the $ L^2 $ estimate, we assume
\begin{equation*}
\left\|W_{\v}^{\r}\right\|^2_{|\a|-1}\leq C\r^{2(2m\added{+1}-(|\a|-1))}= C\r^{2(2m+2-|\a|)},  
\end{equation*}
then, \eqref{ho} becomes
\begin{eqnarray}\label{ho37}
&\quad&\frac{d}{dt}\left( \sum_{\nu=e,i}\big\langle A_\nu^0(n^\r_{\nu})\pa^\a_xW_{\nu}^{\r},\pa^\a_xW_{\nu}^{\r}\big\rangle+\|\nabla \pa_x^\a\Phi^\eps\|^2\right)\nonumber\\
&\leq& C\sum_{\nu=e,i}\|W_\nu^\eps\|_{|\a|}^2+\|\nabla \Phi^\eps\|_{|\a|}^2+C\r^{2(2m\added{+1}-|\a|)}+C\eps^{4m\added{+2}}\nonumber\\.
&\added{\leq}& \added{C\sum_{\nu=e,i}\|W_\nu^\eps\|_{|\a|}^2+\|\nabla \Phi^\eps\|_{|\a|}^2+C\r^{2(2m\added{+1}-|\a|)}.}
\end{eqnarray}
\added{Since $A_\nu^0(n_\nu^\eps)$ is positive definite, consequently, $\big\langle A_\nu^0(n^\r_{\nu})\pa^\a_xW_{\nu}^{\r},\pa^\a_xW_{\nu}^{\r}\big\rangle$ is equivalent to $\|\pa^\a_xW_{\nu}^{\r}\|^2$. Summing \eqref{ho37} up to the $\a$-th order derivative and using the Gronwall inequality, we have}
\begin{equation}\label{final}
\added{\|\nabla\Phi^\eps\|_{|\a|}^2+\left\|W_{\v}^{\r}\right\|^2_{|\a|}\leq  C\r^{2(2m\added{+1}-|\a|)}, }
\end{equation}
\added{so that the induction argument is complete. Combining \eqref{final} for all $1\leq |\a|\leq s$ and noticing \eqref{F3.78}, we have}
\deleted{Together with Lemma \ref{lem1}, we have}
\begin{equation}\label{final2}
\sup_{0\leq t\leq T^{\r,1}}\left\|W_{\added{\nu}}^{\r}(t)\right\|^2_s\leq C\r^{2(2m\added{+1}-s)}, \added{\quad \nu=e,i.}
\end{equation}

It suffices to prove $T_1^{\r,1}\geq T_1^e$, i.e. $T_2^{\r,1}=T_1^e$. By the definitions of  $T_1^e,\,T_1^{\r,1},\,T_2^{\r,1}$ and $T^{\r,1}$, we have $T^{\r,1}\leq T_2^{\r,1}\leq T_1^e$. According to \replaced{the argument of the uniform boundedness of $W_\nu^\eps$}{\eqref{3.6}}, we may replace $T^{\r,1}$ by $T_*^{\r,1}\in(0,\,T_1^e]$ such that $\left[0,\,T_*^{\r,1}\right]$ is the maximum time interval on which $W_{\added{\nu}}^\r$ exists and satisfies \deleted{\eqref{3.6}, i.e.}
\begin{equation*}
\|W_{\added{\nu}}^{\r}\added{(t)}\|_s^2\leq C,\quad\forall\,t\in\left[0,\,T_*^{\r,1}\right],
\end{equation*}
for some constant $C>0$. \added{By \eqref{final2}, it is obvious that} 
\begin{equation*}
\left\|W_{\added{\nu}}^{\r}\left(T_*^{\r,1}\right)\right\|^2_s\leq C\r^{2(2m\added{+1}-s)}.
\end{equation*}
We want to prove $T_*^{\r,1}=T_1^e$. If $T_*^{\r,1}< T_1^e$, we apply the \replaced{theories}{theorem} of Kato for the local existence of smooth solutions with initial data $W_{\added{\nu}}^{\r}\left(T_*^{\r,1}\right)$. Consequently, there exist\added{s} $T_{\r,1}>T_*^{\r,1}$ \added{such that $W_{\added{\nu}}^{\r}\in C\left(\left[T_*^{\eps,1},\,T_{\r,1}\right];H^s\left(\mathbb{R}^d\right)\right)$}\deleted{and a smooth solution $W_{\added{\nu}}^{\r}\in C\left(\left[0,\,T_{\r,1}\right];H^s\left(\mathbb{R}^3\right)\right)$ of \eqref{jh}-\eqref{jhcz1}.} When $2m>s$ and $\r$ is sufficiently small, we always
have $\r^{2(2m\added{+1}-s)}<C$  \deleted{for all fixed constant $C > 0$}. Since the function $t\rightarrow\left\|W_{\added{\nu}}^{\r}(t)\right\|_s$ is continuous on $\left[T_*^{\r,1},\,T_{\r,1}\right]$, there exists $T'_{\r,1}\in\left(T_*^{\r,1},\,T_{\r,1}\right]$ such that
\begin{equation*}
\left\|W_{\added{\nu}}^{\r}(t)\right\|_s^2\leq C,\quad t\in\left[0,\,T'_{\r,1}\right].
\end{equation*}
This is contradictory to the maximality of $T_*^{\r,1}$. Thus, we have proved $T_*^{\r,1}= T_1^e$, which implies that $T_1^{\r,1}\geq T_1^e$.\hfill$\square$
\section{Proof of Theorem \ref{i2.1}}
\subsection{Energy estimates.}
In this section, \added{we study the infinity-ion mass limit of the bipolar Euler-Poisson equations.} We continue to use  $ \left(n_\nu^\r,\,u_\nu^\r,\,\phi^\r\right) $ and $ \left(n_\nu^j,\,u_\nu^j,\,\phi^j\right)_{ j\geq 0} $ to replace $ \left(n_\nu^{1,\frac{1}{\r}},\,u_\nu^{1,\frac{1}{\r}},\,\phi^{1,\frac{1}{\r}}\right) $ and $\left(n_\nu^{i,j},\,u_\nu^{i,j},\,\phi^{i,j}\right)_{ j\geq 0} $. \added{All the corresponding notations are in accordance with those in Subsection 2.3.} The exact solution $\left(n_\v^\r,\,u_\v^\r,\,\phi^\r\right)$
is defined \replaced{on the}{in} time interval $\left[0,\,T_1^{1,\fr} \right]$ and the approximate
solution $\left(n^m_{\v,\r},\,u^m_{\v,\r},\,\phi^m_\r\right)$  \replaced{on the}{in} time
interval $\left[0,\,T_1^i \right]$, with  \added{$T_1^i$ independent of $\eps$}\deleted{$T_1^{1,\fr}>0$ and $T_1^i>0$}. Let
$$      T_2^{1,\fr}=\min\left(T_1^{1,\fr},\,T_1^i\right)>0,  $$
then the exact solution and the approximate solution are both defined \replaced{on the}{in} time interval $\left[0,\,T_2^{1,\fr} \right]$. \replaced{On}{In} this time interval, we denote
\begin{equation*}
\left(N_\v^{\eps},\,U_\v^{\eps},\Phi^{\eps}\right)
\triangleq\left(n_\v^\r-n_{\v,\r}^{m},\,u_\v^\r-u_{\v,\r}^{m},\phi^\r-\phi^m_\r\right),\quad \v=e,\,i.
\end{equation*}
\deleted{For simplicity, we denote
$\left(N_\v^{1,\fr},\,U_\v^{1,\fr},\Phi^{1,\fr},R^{1,\fr,m}\right)$ by $\left(N_\v^{\r},\,U_\v^{\r},\Phi^{\r},R^{\r,m}\right)$ in this section.} \added{For simplicity, set $\lambda=1$.} It is easy to check that \deleted{ the variable}
$\left(N_\v^\r,\,U_\v^\r\right)$ satisfy
\begin{equation}
\label{4E-M3}
\left\{\aligned
&\p N_{e}^{\r}+\big(U_{e}^{\r}+u_{e,\r}^{m}\big)\D N_{e}^{\r}
+\big(N_{e}^{\r}+n_{e,\r}^{m}\big)\div U_{e}^{\r}
=-\big(N_{e}^{\r}\div u_{e,\r}^{m}+U_{e}^{\r}\D n_{e,\r}^{m}\big)-R_{n_e}^{\r,m}, \\
& \p U_{e}^{\r}+ \big(\big(U_{e}^{\r}+u_{e,\r}^{m}\big)\cdot\D\big) U_e^\r
+ h'_{e}(N_e^\r+n_{e,\r}^{m})\D N_e^\r\\
&\hspace{1.5cm}=- \big(U_{e}^{\r}\cdot\D\big) u_{e,\r}^m
- \big(h'_{e}(N_e^\r+n_{e,\r}^m)-h'_{e}(n_{e,\r}^m)\big)\D n_{e,\r}^m  + \nabla\Phi^\r- R_{u_e}^{\r,m},\\
&\p N_{i}^{\r}+\big(U_{i}^{\r}+u_{i,\r}^{m}\big)\D N_{i}^{\r}
+\big(N_{i}^{\r}+n_{i,\r}^{m}\big)\div U_{i}^{\r}
=-\big(N_{i}^{\r}\div u_{i,\r}^{m}+U_{i}^{\r}\D n_{i,\r}^{m}\big)-R_{n_i}^{\r,m}, \\
&\frac{1}{\r}\p U_{i}^{\r}+\frac{1}{\r}\big(\big(U_{i}^{\r}+u_{i,\r}^{m}\big)\cdot\D\big) U_i^\r
+\r h'_{i}(N_i^\r+n_{i,\r}^{m})\D N_i^\r\\
&\hspace{1.5cm}=-\frac{1}{\r}\big(U_{i}^{\r}\cdot\D\big) u_{i,\r}^m
-\r\big(h'_{i}(N_i^\r+n_{i,\r}^m)-h'_{i}(n_{i,\r}^m)\big)\D n_{i,\r}^m -\r\nabla\Phi^\r-\dfrac{1}{\r} R_{u_i}^{\r,m},\\
&\big.\big(N_\v^\r,\,U_\v^\r\big)\big|_{t=0} =\big(n_{\v,0}^\r-n_{\v,\r}^m(0,\,\cdot),
\,u_{\v,0}^\r-u_{\v,\r}^m(0,\,\cdot)\big),
\endaligned
\right.
\end{equation}
with the Poisson equation for  $ \Phi^\r $
\begin{equation*}
    -\Delta \Phi^\r=N_{i}^{\r}-N_{e}^{\r}, \quad \lim_{|x|\rightarrow+\infty}\Phi^\r(x)= 0.
\end{equation*}
Set
$$  W_{ e}^{\r}=\left(\begin{array}{c}
     N_{e}^{\r} \\
       U_{e}^{\r} \\
   \end{array}
 \right),\quad
W_{i}^{\r}=\left(\begin{array}{c}
     N_{i}^{\r} \\
     \frac{1}{\r}U_{i}^{\r} \\
   \end{array}
 \right),
$$
\begin{eqnarray*}
H_{e,\r}^{1}&=&
\left(
   \begin{array}{c}
   N_{e}^{\r}\div u_{e,\r}^{m}+U_{e}^{\r}\D n_{e,\r}^{m} \\
 \left(U_{e}^{\r}\cdot\D\right) u_{e,\r}^{m}+
 \left(h'_{e}(N_e^\r+n_{e,\r}^{m})-h'_{e}(n_{e,\r}^{m})\right)\D n_{e,\r}^{m}\\
   \end{array}
 \right),
\\
H_{i,\r}^{1}&=&
\left(
   \begin{array}{c}
   N_{i}^{\r}\div u_{i,\r}^{m}+U_{i}^{\r}\D n_{i,\r}^{m} \\
\frac{1}{\r}\left(U_{i}^{\r}\cdot\D\right) u_{i,\r}^{m}+\r\left(h'_{i}(N_i^\r+n_{i,\r}^{m})-h'_{i}(n_{i,\r}^{m})\right)\D n_{i,\r}^{m} \\
   \end{array}
 \right),
\end{eqnarray*}
$$
H_{e,\r}^{2}=
\left(
   \begin{array}{c}
0 \\
\nabla\Phi^\r\\
   \end{array}
 \right),
\quad
H_{i,\r}^{2}=
\left(
   \begin{array}{c}
0 \\
-\r\nabla\Phi^\r\\
   \end{array}
 \right),
$$$$
R_{e}^{\r}=
\left(
   \begin{array}{c}
    R_{n_e}^{\r} \\
R_{u_e}^{\r} \\
   \end{array}
 \right) ,
 R_{i}^{\r}=
\left(
   \begin{array}{c}
    R_{n_i}^{\r} \\
\fr R_{u_i}^{\r} \\
   \end{array}
 \right), $$
and for \replaced{$1\leq j\leq d$}{$j=1,2,3$} \added{and $u_\nu^\eps=(u_{\nu,1}^\eps, \cdots, u_{\nu,d}^\eps)$,}
\begin{eqnarray*}
A_e^j\left(n_{e}^{\r},\,u_{e}^{\r}\right)&=&
\left(
  \begin{array}{cc}
u_{e,j}^{\r} &  n_{e}^{\r}e_j^{\top} \\
h'_e\left(n_{e}^{\r}\right)e_j& u_{e,j}^{\r}\mbox{\bf I}_d \\
  \end{array}
\right),\,\\
A_i^j\left(n_{i}^{\r},\,u_{i}^{\r}\right)&=&
\left(
  \begin{array}{cc}
u_{i,j}^{\r} & {\r}n_{i}^{\r}e_j^{\top} \\
{\r}h'_i\left(n_{i}^{\r}\right)e_j&u_{i,j}^{\r}\mbox{\bf I}_d \\
  \end{array}
\right),
\end{eqnarray*}
where $(e_1,\cdots\,e_d)$ is the canonical basis of $\mathbb{R}^d$ and $\mbox{\bf I}_d$
is the $d\times d$ unit matrix. Thus
the  equations \deleted{of} \eqref{4E-M3} can be written as
\begin{equation}\label{5jh}
 \p W_{\v}^{\r}+\sum_{j=1}^d A_\v^j\left(n_{\v}^{\r},\,u_{\v}^{\r}\right)\partial_{x_j}W_{\v}^{\r}
=-H_{\v,\r}^{1} +H_{\v,\r}^{2}-R_{\v}^{\r},\quad \v=e,\,i,
 \end{equation}
with the initial data \deleted{is}
\begin{equation}
\label{jhcz5}
t=0\,:\quad W_{ \v}^{\r}=W_{ \v,0}^{\r},\quad \v=e,\,i,
\end{equation}
where
\begin{eqnarray*}
  W_{ e,0}^{\r}&=& \left(\begin{array}{c}
     N_{e}^{\r}(0,\,\cdot) \\
     {\r}U_{e}^{\r}(0,\,\cdot) \\
   \end{array}\right)
= \left(\begin{array}{c}
     n_{e,0}^\r-n_{e,\r}^m(0,\,\cdot) \\
     u_{e,0}^\r-u_{e,\r}^m(0,\,\cdot)  \\
   \end{array}\right),\\
   W_{ i,0}^{\r}&=& \left(\begin{array}{c}
     N_{i}^{\r}(0,\,\cdot) \\
     U_{i}^{\r}(0,\,\cdot) \\
   \end{array}\right)
= \left(\begin{array}{c}
     n_{i,0}^\r-n_{i,\r}^m(0,\,\cdot) \\
    \frac{1}{\r}\left(u_{i,0}^\r-u_{i,\r}^m(0,\,\cdot)\right)  \\
   \end{array}\right).
\end{eqnarray*}
System \eqref{5jh}-\eqref{jhcz5} for $W_{ \v}^{\r}$ is symmetrizable hyperbolic when $n_{\v}^{\r}>0$, \replaced{which is ensured by Proposition 2.1.}{. Indeed, since the density $n^0$ of the leading profile satisfies}
\replaced{Indeed}{With this}, let
$$   A_e^0\left(n_{e}^{\r}\right)=
\left(\begin{array}{cc}
h'_e\left(n_{e}^{\r}\right) & 0 \\
0& n_{e}^{\r}\mbox{\bf I}_d \\
  \end{array}\right), \quad
  A_i^0\left(n_{i}^{\r}\right)=
\left(\begin{array}{cc}
h'_i\left(n_{i}^{\r}\right) & 0 \\
0&n_{i}^{\r}\mbox{\bf I}_d \\
  \end{array}\right)\added{,}
$$
\replaced{thus}{and} for \replaced{$1\leq j\leq d$}{$j=1,2,3$},
\begin{eqnarray*}
\tilde{A}_e^j\left(n_{e}^{\r},\,u_{e}^{\r}\right)
&=&A_e^0\left(n_{e}^{\r}\right)A_e^j\left(n_{e}^{\r},\,u_{e}^{\r}\right)
= \left(
  \begin{array}{cc}
h'_e\left(n_{e}^{\r}\right)u_{e,j}^{\r} & p'_e\left(n_{e}^{\r}\right) e_j^{\top} \\
 p'_e\left(n_{e}^{\r}\right)e_j&  n_{e}^{\r}u_{e,j}^{\r}\mbox{\bf I}_d \\
  \end{array}
\right),\\
\tilde{A}_i^j\left(n_{i}^{\r},\,u_{i}^{\r}\right)
&=&A_i^0\left(n_{i}^{\r}\right)A_i^j\left(n_{i}^{\r},\,u_{i}^{\r}\right)
=\left(
  \begin{array}{cc}
h'_i\left(n_{i}^{\r}\right)u_{i,j}^{\r} & {\r}p'_i\left(n_{i}^{\r}\right) e_j^{\top} \\
 {\r}p'_i\left(n_{i}^{\r}\right)e_j&n_{i}^{\r}u_{i,j}^{\r}\mbox{\bf I}_d \\
  \end{array}
\right),
\end{eqnarray*}
then for $n_{\v}^{\r}>0$, $A_\v^0$ is positively definite and $\tilde{A}_\v^j$ is
symmetric for all \replaced{$1\leq j\leq d$}{$1\leq j\leq 3$}. Thus, the \replaced{theories}{theorem} of Kato for the local
existence of smooth solutions can also be applied to \eqref{5jh}-\eqref{jhcz5}. \added{Following the similar procedures as those in the section of zero-electron mass limit, we can obtain that there is a time $T^{1,\fr}>0$, such that for $\eps<1$ and $\forall\, t\in[0,T^{1,\fr}]$,
\[
\|W_\nu^\eps(t)\|_{W^{1,\infty}(\R^d)}+\left\|\left(n^m_{\v,\r},\,u^m_{\v,\r},\,\phi^m_\r\right)(t)\right\|_{W^{1,\,\infty}(\R^d)}+\|(n_\nu^\eps, u_\nu^\eps,\phi^\eps )(t)\|_{W^{1,\,\infty}(\R^d)}\leq C.
\]
In order to prove $T_1^{1,\fr}\geq T_1^i$, we need to show that there exists a constant
$\mu > 0$ such that
\begin{equation*}
\sup_{0\leq t\leq T^\r}\|W_{\added{\nu}}^{\r}\added{(t)}\|_s \leq C\r^\mu.
\end{equation*}
}
\subsubsection{\bf $L^2$-estimates.}
In what follows, we always assume that the conditions of Theorem \ref{i2.1} hold.

\begin{lemma}
\label{5lem1}
For all $t\in\left[0,\,T^{1,\fr}\right]$ and sufficiently small $\r>0$, we have
\begin{eqnarray}\label{F3.782}
&\,&\frac{d}{dt}\left(\sum\limits_{\v=e,\,i}\big\langle A_{\v}^0\left(n_{\v}^{\r}\right)W_{\v}^{\r},\,W_{\v}^{\r}
\big\rangle+\|\nabla\Phi^\r\|^2\right)\nonumber\\
&\leq& C\sum\limits_{\v=e,\,i}\left\|W_{\v}^{\r}\right\|^2+\|\nabla\Phi^\r\|^2+C{\r}^{4m\added{+2}}.
\end{eqnarray}
\end{lemma}

\noindent \textbf{Proof. }\underline{Step}1:
Taking the inner product of the \replaced{equations for ions}{ion equations} in \eqref{5jh} with $2A_i^0\left(n_{i}^{\r}\right)W_{i}^{\r}$ in \replaced{$L^2(\mathbb{R}^d)$}{$L^2(\mathbb{R}^3)$}, we obtain the following energy equality for $W_{i}^{\r}$
\begin{align}
\label{5F3.8}
    \frac{d}{dt}\big\langle A_i^0(n^\r_{i})W_{i}^{\r},W_{i}^{\r}\big\rangle
=& \big\langle {\mathrm{div}} A_i(n^\r_{i},u^\r_{i})W_{i}^{\r},W_{i}^{\r}\big\rangle
-2\big\langle A_i^0\left(n_{i}^{\r}\right)W_{i}^{\r},\,H_{i,\r}^{1}\big\rangle
\nonumber\\
&+2\big\langle A_i^0\left(n_{i}^{\r}\right)W_{i}^{\r},\,H_{i,\r}^{2}\big\rangle
-2\big\langle A_i^0\left(n_{i}^{\r}\right)W_{i}^{\r},\,R_{i}^{\r}\big\rangle,
\end{align}
where
\begin{equation}\label{divA_i2}\div A_i\left(n_{i}^{\r},\,u_{i}^{\r}\right)=\p A_i^0\left(n_{i}^{\r}\right)+\sum_{j=1}^d \partial_{x_j} \tilde{A}_i^j\left(n_{i}^{\r},\,u_{i}^{\r}\right).\end{equation}
Now we deal with each term on the right hand side of (\ref{5F3.8}). \added{In the later proof, we will frequently use the fact that for sufficiently small $\eps$, $1<\dfrac{1}{\r}$}. First, from the mass conservation law
 $ \pt n^\r_{i} = - {\mathrm{div}}(n^\r_{i} u^\r_{i}) $,  we have
\added{\[
	\|\pt A_i^0(n_e^\eps)\|_\infty\leq C\|\pt n^\r_{i}\|_\infty\leq C\|{\mathrm{div}}(n^\r_{i} u^\r_{i})\|_{s-1}\leq C, 
	\]which implies}
\begin{eqnarray}
\label{5F3.9}
  \big\langle\pt A_i^0(n^\r_{i})W_{i}^{\r},W_{i}^{\r}\big\rangle&\leq& C\|W_{i}^{\r}\|^2,
\end{eqnarray}
and in view of the expression of  $ \tilde{A}_i^j(W_{i}^{\r}) $,  we obtain
\begin{eqnarray*}
  \big<\pa_{x_j}\tilde{A}_i^j(n_{i}^{\r},\,u_{i}^{\r})W_{i}^{\r},W_{i}^{\r}\big>
&\!\!=\!\!& \big<\pa_{x_j}(h'_{i}(n^\r_{i})u^\r_{i,j})N_{i}^{\r},N_{i}^{\r}\big>
+2\big<N_{i}^{\r}\pa_{x_j}(p'_{i}(n^\r_{i})e_j),{U_{i}^{\r}}\big>\nonumber\\[2mm]
&\!\!\,\!\!& +\, \frac{1}{\r^2}\big<\pa_{x_j}(n^\r_{i} u^\r_{i,j}){U_{i}^{\r}},{U_{i}^{\r}}\big>,
\end{eqnarray*}
in which
\[  \big<\pa_{x_j}(h'_{i}(n^\r_{i})u^\r_{i,j})N_{i}^{\r},N_{i}^{\r}\big>
+ \frac{1}{\r^2}\big<\pa_{x_j}(n^\r_{i} u^\r_{i,j}){U_{i}^{\r}},{U_{i}^{\r}}\big>
\leq C\|W_{i}^{\r}\|^2,  \]
and
\begin{eqnarray*}
2\sum_{j=1}^d\big<N_{i}^{\r}\pa_{x_j}(p'_{i}(n^\r_{i})e_j),{U_{i}^{\r}}\big>
=\added{2}\big\langle N_{i}^{\r}\D p'(n^\r_{i}),{U_{i}^{\r}}\big\rangle
&\leq& C\|U_{i}^{\r}\|^2+C\|N_{i}^{\r}\|^2\nonumber\\
&\leq& \frac{C}{\r^2}\|U_{i}^{\r}\|^2+C\|N_{i}^{\r}\|^2\nonumber\\
&\leq& C\|W_{i}^{\r}\|^2,
\end{eqnarray*}
therefore,
\begin{equation}
\label{5F3.10}
 \sum_{j=1}^d\big\langle\pa_{x_j}\tilde{A}_i^j(n_{i}^{\r},\,u_{i}^{\r})W_{i}^{\r},W_{i}^{\r}\big\rangle
\leq C\|W_{i}^{\r}\|^2.
\end{equation}
It follows from   \eqref{5F3.9} and \eqref{5F3.10} that
\begin{equation}
\label{5F3.11}
   \big\langle {\mathrm{div}} A_i(n^\r_{i},u^\r_{i})W_{i}^{\r},W_{i}^{\r}\big\rangle
\leq C\|W_{i}^{\r}\|^2.
\end{equation}
For the remaining terms without $H_{e,\eps}^2$ in the right hand side of (\ref{5F3.8}), we have
\begin{eqnarray*}
&\,&-2\big\langle A_i^0(n^\r_{i})W_{i}^{\r},H_{i,\r}^{1}\big\rangle
-2\big\langle A_i^0\left(n_{i}^{\r}\right)W_{i}^{\r},\,R_{i}^{\r}\big\rangle
\nonumber\\
&=& - 2\big\langle h_i'(n^\r_{i})
\big(N_{i}^{\r}{\mathrm{div}} u_{i,\r}^m + {U_{i}^{\r}}\D n_{i,\r}^m
+R_{n_i}^{\r,m}\big),N_{i}^{\r}\big\rangle
 \nonumber\\
&&-\, 2\Big\langle n^\r_{i}\Big[\dfrac{1}{\r^2}\big({U_{i}^{\r}}\cdot\D\big) u_{i,\r}^m+\left(h'_{i}(N_i^\r+n_{i,\r}^{m})-h'_{i}(n_{i,\r}^{m})\right)\D n_{i,\r}^{m}
+ \dfrac{1}{\eps^2}R_{u_i}^{\r,m}\Big],{U_{i}^{\r}}\Big\rangle\nonumber\\
&=& -\,2\big\langle h_i'(n^\r_{i})N_{i}^{\r}{\mathrm{div}} u_{i,\r}^m,N_{i}^{\r}\big\rangle
- \frac{2}{\r^2}\big\langle n^\r_{i}\big({U_{i}^{\r}}\cdot\D\big) u_{i,\r}^m,{U_{i}^{\r}}\big\rangle
- 2\big\langle h_i'(n^\r_{i})R_{n_i}^{\r,m},N_{i}^{\r}\big\rangle\nonumber\\
&\,&- \dfrac{2}{\eps^2}\big\langle n^\r_{i} R_{u_i}^{\r,m},{U_{i}^{\r}}\big\rangle
-\,2\big\langle h_i'(n^\r_{i})N_{i}^{\r}\D n_{i,\r}^m,{U_{i}^{\r}}\big\rangle
- 2 \big\langle n^\r_{i}\big(h_i'(n^\r_{i})-h_i'(n_{i,\r}^m)\big)
\D n_{i,\r}^m,{U_{i}^{\r}}\big\rangle,
\end{eqnarray*}
in which \added{by the Taylor's formula,}
\begin{eqnarray*}
 -\,2\big\langle h_i'(n^\r_{i})N_{i}^{\r}{\mathrm{div}} u_{i,\r}^m,N_{i}^{\r}\big\rangle
- \frac{2}{\r^2}\big\langle n^\r_{i}\big({U_{i}^{\r}}\cdot\D\big) u_{i,\r}^m,{U_{i}^{\r}}\big\rangle
&\leq & C\|W_{i}^{\r}\|^2,\\
 -2\big\langle h_i'(n^\r_{i})N_{i}^{\r}\D n_{i,\r}^m, U_i^\eps\big\rangle
-2\big\langle n^\r_{i}\big(h_i'(n_{i,\r}^m)-h_i'(n^\r_{i})\big)\D n_{i,\r}^m,{U_{i}^{\r}}\big\rangle
&\leq& C\|W_{i}^{\r}\|^2,
\end{eqnarray*}
and
\[
- 2\big\langle h_i'(n^\r_{i})R_{n_i}^{\r,m},N_{i}^{\r}\big\rangle
- \dfrac{2}{\eps^2}\big\langle n^\r_{i} R_{u_i}^{\r,m},{U_{i}^{\r}}\big\rangle
\leq  C\|W_{i}^{\r}\|^2
+ C\Big( \|R_{n_i}^{\r,m}\|^2 + \added{\dfrac{1}{\eps^2}}\|R_{u_i}^{\r,m}\|^2\Big).\\
\]
As for the term containing $H_{i,\r}^{2}$ in \eqref{5F3.8}, a direct calculation gives
\begin{align*}\label{5F3.15}
2\big\langle  A_i^0\left(n_{i}^{\r}\right)W_{i}^{\r},\,H_{i,\r}^{2}\big\rangle
=-2\big\langle n_{i}^{\r}U_{i}^{\r},\,\nabla\Phi^\r\big\rangle.
\end{align*}
Finally, using \eqref{5F3.8}, \eqref{5F3.11} and the four estimates above yield
\begin{eqnarray}
\label{5F3.16}
 &\,&\frac{d}{dt} \big\langle A_i^0(n^\r_{i})W_{i}^{\r},W_{i}^{\r}\big\rangle \nonumber\\
  &\leq& C\|W_{i}^{\r}\|^2-2\big\langle n_{i}^{\r}U_{i}^{\r},\,\nabla\Phi^\r\big\rangle
+ C\Big( \|R_{n_i}^{\r,m}\|^2 + \added{\dfrac{1}{\eps^2}}\|R_{u_i}^{\r,m}\|^2\Big).
\end{eqnarray}

\underline{Step}2: Similar to what we have done in the previous section, taking the inner product of the \replaced{equations for electrons}{electrons equations} in \eqref{5jh} with $2A_e^0\left(n_{e}^{\r}\right)W_{e}^{\r}$ in \replaced{$L^2(\mathbb{R}^d)$}{$L^2(\mathbb{R}^3)$}, we have
\begin{align*}
\frac{d}{dt}\big\langle A_e^0\left(n_{e}^{\r}\right)W_{e}^{\r},\,W_{e}^{\r}\big\rangle
=&-2\big\langle A_e^0\left(n_{e}^{\r}\right)W_{e}^{\r},\,H_{e,\r}^{1}\big\rangle
+2\big\langle A_e^0\left(n_{e}^{\r}\right)W_{e}^{\r},\,H_{e,\r}^{2}\big\rangle\nonumber\\
&-2\big\langle A_e^0\left(n_{e}^{\r}\right)W_{e}^{\r},\,R_{e}^{\r}\big\rangle
+\big\langle\div A_e\left(n_{e}^{\r},\,u_{e}^{\r}\right)W_{e}^{\r},\,W_{e}^{\r}\big\rangle,
\end{align*}
where
\begin{equation}\label{divA_e2}\div A_e\left(n_{e}^{\r},\,u_{e}^{\r}\right)=\p A_e^0\left(n_{e}^{\r}\right)+\sum_{j=1}^d \partial_{x_j} \tilde{A}_e^j\left(n_{e}^{\r},\,u_{e}^{\r}\right). \end{equation}
The estimates are all the same as we did \added{for the equations for ions} in the zero-electron mass limit, since both of them do not \replaced{contain}{have} the parameters $\r$ and \replaced{are}{is} only different in notations, we omit the proof. Indeed, we have
\begin{equation}\label{imlm2}
\frac{d}{dt}\big\langle A_e^0\left(n_{e}^{\r}\right)W_{e}^{\r},\,W_{e}^{\r}\big\rangle
\leq2\big\langle n_{e}^{\r}U_{e}^{\r},\,\nabla\Phi^\r\big\rangle+C\left\|W_{e}^{\r}\right\|^2+C{\r}^{4m\added{+4}} .
\end{equation}

\underline{Step}3: Summing \eqref{5F3.16} and \eqref{imlm2} \deleted{for all $|\a|\leq s$}, following \replaced{a similar}{ the same} procedure as the $L^2$-estimate in the previous section, we obtain \eqref{F3.782}. \hfill$\square$
\subsubsection{\bf Higher order estimates.}
Let \replaced{$\a\in\mathbb{N}^d$}{$\a\in\mathbb{N}^3$} with $1\leq|\a|\leq s$. Applying $\partial^\a_x$ to \eqref{5jh}, we get
\begin{equation}\label{5jh2}
\p \partial^\a_xW_{\v}^{\r}+\sum_{j=1}^d A_{\v}^j\left(n_{\v}^{\r},\,u_{\v}^{\r}\right)\partial_{x_j}\partial^\a_xW_{\v}^{\r}
=-\partial^\a_x\left(H_{\v,\r}^{1}-H_{\v,\r}^{2}+R_{\v}^{\r}\right)
+J_{\v,\r}^{\a},\quad \v=e,\,i,
\end{equation}
where$$J_{\v,\r}^{\a}=\sum_{j=1}^d \left(A_{\v}^j\left(n_{\v}^{\r},\,u_{\v}^{\r}\right)\partial_{x_j}\partial^\a_xW_{\v}^{\r}
-\partial^\a_x\left(A_{\v}^j\left(n_{\v}^{\r},\,u_{\v}^{\r}\right)\partial_{x_j}W_{\v}^{\r}\right)\right).$$

\begin{lemma}
\label{5lem3}
   For all $t\in\left[0,\,T^{1,\fr}\right]$ and sufficiently small $\r>0$, we have
   \begin{eqnarray}\label{ho2}
&\,&\frac{d}{dt}\left( \sum_{\nu=e,i}\big\langle A_\nu^0(n^\r_{\nu})\pa^\a_xW_{\nu}^{\r},\pa^\a_xW_{\nu}^{\r}\big\rangle+\|\nabla \pa_x^\a\Phi^\eps\|^2\right)\nonumber\\
&\leq& C\sum_{\nu=e,i}\|W_\nu^\eps\|_{|\a|}^2+\|\nabla \Phi^\eps\|_{|\a|}^2+C\eps^{4m\added{+2}}.
\end{eqnarray}

\end{lemma}

\noindent\textbf{Proof.} \underline{Step}1:Taking the inner product of the \replaced{equations for ions}{ion equations} in \eqref{5jh2} with
$2A_i^0\left(n_{i}^{\r}\right)\partial^\a_xW_{i}^{\r}$ in \replaced{$L^2(\mathbb{R}^d)$}{$L^2(\mathbb{R}^3)$}
yields the following energy equality for $\partial^\a_xW_{i}^{\r}$
\begin{eqnarray}
\label{5F3.23}
&\quad&\frac{d}{dt}\big\langle A_i^0(n^\r_{i})\pa^\a_xW_{i}^{\r},\pa^\a_xW_{i}^{\r}\big\rangle\nonumber\\
&=& \big\langle {\mathrm{div}} A_i(n^\r_{i},u^\r_{i})\pa^\a_xW_{i}^{\r},\pa^\a_xW_{i}^{\r}\big\rangle-2\big\langle A_i^0\left(n_{i}^{\r}\right)\pa^\a_xW_{i}^{\r},\,\pa^\a_xH_{i,\r}^{1}+\pa^\a_xR_{i}^{\r}\big\rangle\nonumber\\
&\quad&+2\big\langle A_i^0\left(n_{i}^{\r}\right)\pa^\a_xW_{i}^{\r},\,\pa^\a_xH_{i,\r}^{2}\big\rangle+2\left\langle A_i^0\left(n_{i}^{\r}\right)\partial^\a_xW_{i}^{\r},\,J_{i,\r}^{\a}\right\rangle,
\end{eqnarray}
\added{where $\dive A_i$ is defined in \eqref{divA_i2}. We will treat the right hand side of the above} \deleted{which are treated} term by term as follows. First, similarly to \eqref{5F3.11}, it is easy to get
\begin{equation}
\label{5F3.24}
\left|\left\langle{\mathrm{div}} A_i(n^\r_{i},u^\r_{i})\pa^\a_xW_{i}^{\r},\pa^\a_xW_{i}^{\r}\right\rangle\right|
\leq C\left\|W_{i}^{\r}\right\|^2_{|\a|}.
\end{equation}
For the terms without $H_{i,\r}^2$ and $J_{i,\r}^\a$ \replaced{on}{in} the right hand side of \eqref{5F3.23}, a straightforward calculation yields
\begin{eqnarray*}
&\,& -2\big\langle A_i^0(n^\r_{i})\pa^\a_xW_{i}^{\r},\pa^\a_xH_{i,\r}^{1}+\pa^\a_xR_{i}^{\r}\big\rangle\nonumber\\
&=& - 2\big< h_i'(n^\r_{i}) \big(\pa^\a_x\big(N_{i}^{\r}\dive u_{i,\r}^m \big)
+\pa^\a_xR_{n_{\added{i}}}^{\r,m}\big),\pa^\a_xN_{i}^{\r}\big>
- 2\big< h_i'(n^\r_{i})
(\pa^\a_xN_{i}^{\r} \D n_{i,\r}^m),\pa^\a_xU_{i}^{\r}\big>\nonumber\\
&\,& - 2\big< h_i'(n^\r_{i})\big(\pa^\a_x\big({U_{i}^{\r}}\D n_{i,\r}^m\big)
-\pa^\a_xU_{i}^{\r} \D n_{i,\r}^m\big),\pa^\a_xN_{i}^{\r}\big>\nonumber\\
&\,& -\frac{2}{\r^2}\big< n^\r_{i}\pa^\a_x\Big[\big({U_{i}^{\r}}\cdot\D\big) u_{i,\r}^m
+ R_{u_{\added{i}}}^{\r,m}\Big],\pa^\a_xU_{i}^{\r}\big>
- 2 \left< n^\r_{i}\pa^\a_x\big(h_i'(n^\r_{i})-h_i'(n_{i,\r}^m)\big)
\D n_{i,\r}^m,\pa^\a_xU_{i}^{\r}\right>\nonumber\\
&\,&- 2\big<n^\r_{i} \big[\pa^\a_x\big((h_i'(n^\r_{i})-h_i'(n_{i,\r}^m))
\D n_{i,\r}^m \big)-\pa^\a_x\big(h_i'(n^\r_{i})-h_i'(n_{i,\r}^m)\big)
\D n_{i,\r}^m\big],\pa^\a_xU_{i}^{\r}\big>,
\end{eqnarray*}
to which applying the Moser-type inequalities yields
\begin{eqnarray*}
- 2\left< h_i'(n^\r_{i}) \big(\pa^\a_x(N_{i}^{\r}\dive u_{i,\r}^m )
+\pa^\a_xR_{n_{\added{i}}}^{\r,m}\big),\pa^\a_xN_{i}^{\r}\right>
&\leq&  C\Big(\|W_{i}^{\r}\|^2_{|\a|}+\|R_{n_{\added{i}}}^{\r,m}\|^2_{|\a|}\Big),\\
\label{5F13.13}- 2\big< h_i'(n^\r_{i})
(\pa^\a_xN_{i}^{\r} \D n_{i,\r}^m),\pa^\a_xU_{i}^{\r}\big>
&\leq& C\|W_{i}^{\r}\|^2_{|\a|},\\
\label{F3.28}
  - 2\left< h_i'(n^\r_{i})\big(\pa^\a_x({U_{i}^{\r}}\D n_{i,\r}^m)
-\pa^\a_xU_{i}^{\r} \D n_{i,\r}^m\big),\pa^\a_xN_{i}^{\r}\right>
&\leq &C\|W_{i}^{\r}\|^2_{|\a|},\\
    -\frac{2}{\r^2}\left< n^\r_{i}\pa^\a_x\Big[\big({U_{i}^{\r}}\cdot\D\big) u_{i,\r}^m
+ R_{u_{\added{i}}}^{\r,m}\Big],\pa^\a_xU_{i}^{\r}\right>
&\leq &C\Big(\|W_{i}^{\r}\|^2_{|\a|}
+ \added{\dfrac{1}{\eps^2}}\|R_{u_{\added{i}}}^{\r,m}\|^2_{|\a|}\Big),\\
\label{5F23.13}
- 2 \left< n^\r_{i}\pa^\a_x\big(h_i'(n^\r_{i})-h_i'(n_{i,\r}^m)\big)
\D n_{i,\r}^m,\pa^\a_xU_{i}^{\r}\right>
&\leq& C\|W_{i}^{\r}\|^2_{|\a|},
\end{eqnarray*}
and \added{further by the Taylor's formula,}
\begin{align*}
-2\left< n^\r_{i} \big(\pa^\a_x(\big(h_i'(n^\r_{i})-h_i'(n_{i,\r}^m)\big)
\D n_{i,\r}^m )-\pa^\a_x\big(h_i'(n^\r_{i})-h_i'(n_{i,\r}^m)\big)
\D n_{i,\r}^m\big),\pa^\a_xU_{i}^{\r}\right>
\leq C\|W_{i}^{\r}\|^2_{|\a|}.
\end{align*}
These estimates imply
\begin{eqnarray}
-2\big\langle A_i^0(n^\r_{i})\pa^\a_xW_{i}^{\r},\pa^\a_xH_{i,\r}^{1}+\pa^\a_xR_{i}^{\r}\big\rangle\leq C\left(\|W_{i}^{\r}\|^2_{|\a|} + \|R_{n_{\added{i}}}^{\r,m}\|^2_{|\a|} + \added{\dfrac{1}{\eps^2}}\|R_{u_{\added{i}}}^{\r,m}\|^2_{|\a|}\right).
\end{eqnarray}
For the term containing  $ J_{\added{i,}\r}^{\a} $  \replaced{on}{in} the right hand side of \eqref{5F3.23}, we have \added{for $1\leq j\leq d$ and $U_i^\eps=(U_{i,1}^\eps,\cdots U_{i,d}^\eps)$,}
\begin{eqnarray*}\label{5jam}
&\;& \left<A_i^0(n^\r_{i})\pa^\alpha_xW_i^\r,\pa^\alpha_x (A_i^j(n^\r_{i},u^\r_{i})\pa_{x_j}W_i^\r)
- A_i^j(n^\r_{i},u^\r_{i})\pa^\alpha_x (\pa_{x_j} W_i^\r)\right> \nonumber\\[2mm]
& = & \left<h_i'(n^\r_{i}) \left(\pa^\alpha_x (u^\r_{i,j} \pa_{x_j} N_i^\r)
- u^\r_{i,j} \pa^\alpha_x \pa_{x_j} N_i^\r\right),\pa^\alpha_x N_i^\r\right>
\nonumber\\[2mm]
& \, & + \;  \frac{1}{\r^2}\left<n^\r_{i} \left(\pa^\alpha_x (u^\r_{i,j} \pa_{x_j} U_i^\eps)
- u^\r_{i,j}\pa^\alpha_x\pa_{x_j} U_i^\eps\right),\pa^\alpha_x U_i^\eps\right> \nonumber\\[2mm]
& \, & + \; \left<n^\r_{i} \left(\pa^\alpha_x(h_i'(n^\r_{i})\pa_{x_j} N_i^\r)
-h_i'(n^\r_{i})\pa^\alpha_x\pa_{x_j} N_i^\r\right),\pa^\alpha_x U_{i,j}^\r\right>\nonumber\\[2mm]
& \, & + \;\left<h_i'(n^\r_{i})\left(\pa^\alpha_x(n^\r_{i} \pa_{x_j} U_{i,j}^\r)
- n^\r_{i} \pa^\alpha_x \pa_{x_j} U_{i,j}^\r \right),\pa^\alpha_x N_i^\r\right>\nonumber\\[2mm]
& \leq & C \|W_i^\r\|^2_{|\a|}
+C\left(\|N_i^\r\|^2_{|\a|-1} + \|U_i^\r\|^2_{|\a|}\right)
+C\left(\|N_i^\r\|^2_{|\a|} + \|U_i^\r\|^2_{|\a|-1}\right),
\end{eqnarray*}
which implies
\begin{equation}
\label{5i220}
\left|\left<A_i^0\left(n_{i}^{\r}\right)\partial^\a_xW_{i}^{\r},\,J_{i,\r}^{\a}\right>\right|\leq C\left\|W_{i}^{\r}\right\|^2_{|\a|}.
\end{equation}
For the term containing $H_{i,\r}^{2}$ in the right hand side of \eqref{5F3.23}, a direct calculation gives
\begin{align}\label{5Fe3.15}
2\left<A_i^0\left(n_{i}^{\r}\right)\pa^\a_xW_{i}^{\r},\,\pa^\a_xH_{i,\r}^{2}\right>
=\added{-}2\left<n_{i}^{\r}\pa^\a_xU_{i}^{\r},\,\nabla\pa^\a_x\Phi^\r\right>.
\end{align}
Therefore, using \eqref{5F3.23}-\eqref{5Fe3.15} yield
\begin{eqnarray}
\label{5F3.31}
\frac{d}{dt}\big\langle A_i^0(n^\r_{i})\pa^\a_xW_{i}^{\r},\pa^\a_xW_{i}^{\r}\big\rangle
\leq \added{-}2\left<n_{i}^{\r}\pa^\a_xU_{i}^{\r},\,\nabla\pa^\a_x\Phi^\r\right>
+C \|W_i^\r\|^2_{|\a|}+C\eps^{4m\added{+2}}.
\end{eqnarray}

\underline{Step}2: Taking the inner product of the \replaced{equations for electrons}{electrons equations} in \eqref{5jh2} with
$2A_e^0\left(n_{e}^{\r}\right)\partial^\a_xW_{e}^{\r}$ in \replaced{$L^2(\mathbb{R}^d)$}{$L^2(\mathbb{R}^3)$}
yields the following energy equality for $\partial^\a_xW_{e}^{\r}$
\begin{align*}
\frac{d}{dt}\big\langle A_e^0\left(n_{e}^{\r}\right)\partial^\a_xW_{e}^{\r},\,\partial^\a_xW_{e}^{\r}\big\rangle
=&-2\big\langle A_e^0\left(n_{e}^{\r}\right)\partial^\a_xW_{e}^{\r},\,\partial^\a_xH_{e,\r}^{1}\big\rangle
+2\big\langle A_e^0\left(n_{e}^{\r}\right)\partial^\a_xW_{e}^{\r},\,\partial^\a_xH_{e,\r}^{2}\big\rangle\nonumber\\
&-2\big\langle A_e^0\left(n_{e}^{\r}\right)\partial^\a_xW_{e}^{\r},\,\partial^\a_xR_{e}^{\r}\big\rangle
+2\big\langle A_e^0\left(n_{e}^{\r}\right)\partial^\a_xW_{e}^{\r},\,J_{e,\r}^{\a}\big\rangle\nonumber\\
&+\big\langle\div A_i\left(n_{e}^{\r},\,u_{e}^{\r}\right)\partial^\a_xW_{e}^{\r},\,\partial^\a_xW_{e}^{\r}\big\rangle,
\end{align*}
\added{where $\dive A_e$ is defined in \eqref{divA_e2}.} The estimates are all the same as we did \added{for the equations for ions }in the zero-electron mass limit, since both of them do not \replaced{contain}{have } the parameters $\r$ and \replaced{are}{is} only different in notations. We omit the proof. Indeed, we have
\begin{equation}
\label{lm31}
   \frac{d}{dt}\big\langle A_e^0\left(n_{e}^{\r}\right)\partial^\a_xW_{e}^{\r},\,\partial^\a_xW_{e}^{\r}\big\rangle
\leq2\big\langle n_{e}^{\r}\partial^\a_xU_{e}^{\r},\,\nabla\partial^\a_x\Phi^\r\big\rangle+C\left\|W_{e}^{\r}\right\|^2_{|\a|}+C{\r}^{4m\added{+4}}.
\end{equation}

\vspace{5mm}

\underline{Step}3: Summing \eqref{5F3.31} and \eqref{lm31} \deleted{for all $|\alpha|\leq s $}, following the same procedure as \added{we did for }the higher order estimates in the previous section, we obtain \eqref{ho2}. \hfill$\square$

\subsection{Proof of Theorem \ref{i2.1}} \added{Summing up \eqref{ho2} for $|\a|\leq s$ and combining \eqref{F3.782}, we have}
 \begin{eqnarray*}
&\,&\added{\frac{d}{dt}\left( \sum_{|\a|\leq s}\sum_{\nu=e,i}\big\langle A_\nu^0(n^\r_{\nu})\pa^\a_xW_{\nu}^{\r},\pa^\a_xW_{\nu}^{\r}\big\rangle+\|\nabla \Phi^\eps\|_s^2\right)}\nonumber\\
&\added{\leq}& \added{C\sum_{\nu=e,i}\|W_\nu^\eps\|_{s}^2+\|\nabla \Phi^\eps\|_{s}^2+C\eps^{4m\added{+2}}.}
\end{eqnarray*}
\added{Since $A_\nu^0(n_\nu^\eps)$ is positive definite, $\sum_{|\a|\leq s}\big\langle A_\nu^0(n^\r_{\nu})\pa^\a_xW_{\nu}^{\r},\pa^\a_xW_{\nu}^{\r}\big\rangle$ is equivalent to $\|W_\nu^\eps\|_s^2$. By applying the Gronwall inequality, we have
\[
\|W_\nu^\eps(t)\|_s^2\leq C\eps^{4m\added{+2}}, \quad \forall\, t\in[0,T^{1,\fr}].
\]
}
The rest of the proof is also based on the continuous method, which is similar as what we did in \added{the section of the} zero-electron mass limit, we omit it here.

\subsection*{\bf Acknowledgments.}
\added{Shuai Xi's research was supported by the Cultivation Project of Young and Innovative Talents in Universities of Shandong Province.}

%

\end{document}